\documentclass[12pt]{amsart}
\usepackage[centertags]{amsmath}

\usepackage{amsfonts}

\usepackage{amssymb}

\title[Cofinal Extensions]{The Diversity of \\ Minimal Cofinal Extensions}
\author{James H. Schmerl}

\date{\today}
\def\into{\longrightarrow}

\def\harp{\hspace{-4pt} \upharpoonright \hspace{-4pt}}

\def\pa{{\sf PA}}

\def\bbone {1 {\hspace{-4.6pt} 1}}
\def\bbzero{0 {\hspace{-4.6pt} 0}}

\newcommand{\KK}{{\mathcal K}}
\newcommand{\NN}{{\mathcal N}}
\newcommand{\MM}{{\mathcal M}}

\newcommand{\ZZ}{{\mathbb Z}}

  \DeclareMathOperator{\ssy}{SSy}

         \DeclareMathOperator{\gcis}{GCIS}

      \DeclareMathOperator{\dom}{dom}
       \DeclareMathOperator{\Bin}{Bin}

 \DeclareMathOperator{\cod}{Cod}

  \DeclareMathOperator{\Th}{Th} 
  
  \DeclareMathOperator{\Def}{Def}

   \DeclareMathOperator{\eq}{Eq}

  \DeclareMathOperator{\Lt}{Lt}

\begin{document}

{\abstract{ Fix a countable nonstandard model $\MM$ of Peano Arithmetic. Even with some rather severe restrictions placed on the types of minimal cofinal extensions  $\NN \succ \MM$ that are allowed, we still find that there are $2^{\aleph_0}$ possible theories of $(\NN,M)$ for such~$\NN$'s.}}

\maketitle

The script letters  $\MM, \NN, \KK$ (possibly adorned) always denote models of Peano  Arithmetic (\pa) having domains $M, N, K$, respectively. The set of parametrically definable subsets of $\MM$ is $\Def(\MM)$. If $J \subseteq M$, then $\cod(\MM / J) = \{A \cap J : A \in \Def(\MM)\}$. A {\em cut} of $\MM$ is a subset $J \subseteq M$ such that $0 \in J \neq M$ and if $a \leq b \in J$, then $a+1 \in J$. The cut $J$ is {\em exponentially closed} if $2^a \in J$ whenever $a \in J$. 

Suppose that $\MM \prec \NN$. Their  {\em {\underline G}reatest {\underline C}ommon {\underline I}nitial {\underline S}egment} is 
$$
\gcis(\MM,\NN) = \{b \in M : {\mbox{ whenever }}  \NN \models a \leq b, {\mbox{ then }}  a \in M\},
$$ 
which is  $M$ if $\NN$ is an end extension of $\MM$ and is a cut otherwise.
If $J$ is a cut of $\MM$, then $\NN$ {\em fills} $J$ if there is $b \in N$ such that whenever $a \in J$ and $c \in M \backslash J$, then $\NN \models a < b < c$.  
The {\em interstructure lattice} is $\Lt(\NN / \MM) = \{\KK  : \MM \preccurlyeq \KK \preccurlyeq \NN\}$, ordered by elementary extension. If $1 \leq n < \omega$, then ${\mathbf n}$ is the lattice that is a chain  of $n$ elements. 

One of the themes of \cite{bases} is the diversity of cofinal extensions, exemplified by   the following theorem.

\bigskip

{\sc Theorem} A:  (\cite[Theorem 7.1]{bases}) {\em If $J$ is an exponentially closed cut of countable $\MM$, then there is a set ${\mathcal C}$ of cofinal elementary extensions of $\MM$ such that$:$

\begin{itemize}

\item[(1)] $|{\mathcal C}| = 2^{\aleph_0};$

\item[(2)] if $\NN \in {\mathcal C}$, then $\gcis(\MM,\NN) = J$,   $\cod(\NN / J) = \cod(\MM / J)$ and $\NN$ does not fill $J;$

\item[$(3)$] if $\NN_1,\NN_2 \in {\mathcal C}$ are distinct, then $\Th(\NN_1,M) \neq \Th(\NN_2,M);$

\item[(4)] $\Lt(\NN / \MM) \cong {\mathbf 3}$ for each $\NN \in {\mathcal C}. \ ($\cite[page 285]{bases}$)$

\end{itemize}}

\bigskip

 It was left open, and specifically asked (\cite[Question 7.5]{bases}), whether the ${\mathbf 3}$ in (4) can be replaced by ${\mathbf 2}$  (so that every  $\NN \in {\mathcal C}$ is a {\em minimal} elementary extension of $\MM$).  
 The purpose of this note is to show, in Theorem~B,  that it can,  demonstrating   the diversity of minimal cofinal extensions.

 We let ${\mathcal L}$ be one of the usual finite languages appropriate for $\pa$, say ${\mathcal L} = \{+,\times,0,1,\leq\}$. 
Let ${\mathcal L}^* = {\mathcal L} \cup \{M\}$, where $M$ is a new unary predicate symbol,  so that ${\mathcal L}^*$ is the language appropriate for structures of the form $(\NN,M)$, where $\MM \prec \NN$. 
The theories in (3) of Theorem~A are ${\mathcal L}^*$-theories. 
 
 \bigskip

{\sc Theorem} B:   {\em  There is an ${\mathcal L}^*$-formula $\varphi(x)$ such that whenever  $J$ is an exponentially closed  cut of some countable $\MM$,  then there is a {\em minimal} extension $\NN \succ \MM$ such that  $\gcis(\MM,\NN) = J$,   $\cod(\NN / J) = \cod(\MM / J)$, $\NN$ does not fill $J$ and $\varphi(x)$ defines $I$ in $(\NN,M)$.}

\bigskip

Of course, since there are $2^{\aleph_0}$ subsets $I \subseteq \omega$, there must be a ${\mathcal C}$ as in Theorem~A but with the ${\mathbf 3}$ in (4) being replaced with ${\mathbf 2}$. 

Incidentally, Theorem~A is an easy consequence of Theorem~B. To see why, let $I \subseteq \omega$ and then let $\NN_I$ be   an $\NN$ as in Theorem~B.  Let $a_I$ generate $\NN_I$ over $\MM$. Then, one easily gets a minimal, cofinal  extension $\NN'_I \succ \NN_I$ such that $a_I \in \gcis(\NN'_I,\NN_I)$, where $b_I$ generates $\NN'_I$ over $\NN_I$ and $a_I$ is in the elementary submodel of $\NN'_I$ generated by $b_I$. Since $N_I$ is uniformly definable in $\NN'_I$, then $I$ is uniformly definable in $\NN'_I$. Then, ${\mathcal C} = \{\NN_I' : I \subseteq \omega\}$ is as required by Theorem~A.

There are four sections following this introduction. The first section contains some preliminary material. The second section introduces and discusses humble extensions. Included in that section is a  result (Corollary~2.12) in 
 the spirit of Theorems~A and~B that has  not appeared elsewhere although it is likely that it is known. The third section concerns some purely finite combinatorics that, when interpreted in \pa, will be applied in  \S4, where  the proof of Theorem~B is completed.

\bigskip

 %
 %

{\section{Preliminaries}} This preliminary section comprises three subsections. The first of these fixes some notation and terminology. The second intends to put Theorem~B into context by presenting some known results with proofs that will be used later on. The third subsection discusses some background on interstructure lattices.

\smallskip

%
%

\subsection{Notation} We begin this section by confirming some  notation for models of \pa. As we already said in the introduction, $\MM, \NN, \ldots$ always denote models of $\pa$ with domains $M,N, \ldots$. The standard model ${\mathbb N} = (\omega, \ldots)$ is assumed to be a substructure of every nonstandard model. We distinguish $\MM \prec  \NN$ from $\MM \preccurlyeq \NN$, the latter allowing that $\MM = \NN$ whereas the former does not. We write $\MM \prec_{\sf cf} \NN$ when $\NN$ is a cofinal extension of $\MM$. Cuts are defined in the introduction. A cut $J$ is {\em multiplicatively closed} if $a^2 \in J$ whenever $a \in J$. 

From the introduction, ${\mathcal L}$ is a language appropriate for $\pa$ and ${\mathcal L}^*$ is ${\mathcal L} \cup \{M\}$. If $\MM$ is a model and $C \subseteq M$, then ${\mathcal L}(C)$ and ${\mathcal L}^*(C)$ are the languages ${\mathcal L}$ and ${\mathcal L}^*$ augmented with constant symbols for elements of~$C$.

If $n < \omega$, then $n = \{0,1,\ldots,n-1\}$. There are two Pigeon Hole Principles that we will use:

\smallskip

\begin{itemize}

\item[]

\begin{itemize}

\item[{\sf PHP}1:] If $a,c < \omega$ and  $f : ca + 1 \into c$, then there is $A \subseteq ca + 1$ such that $|A| = a+1$ and $f$ is constant on~$A$.

\smallskip

\item[{\sf PHP2}:] If $a < \omega$ and $f : a^2+1 \into \omega$, then there is $A \subseteq a^2 + 1$ such that $|A| = a+1$ and $f$ is either constant or one-to-one on~$A$.

\end{itemize}
\end{itemize}
\smallskip
Each of these is formalizable and provable in \pa. The formalized versions will be referred to as {\sf PHP1}(\pa) and {\sf PHP2}(\pa).

%
%

\subsection{Some Context}To help put Theorem~B into perspective, we mention several results that are probably well known. 

\bigskip

{\sc Theorem 1.1}:  {\em Suppose that $J$ is a cut of countable $\MM$. The following are equivalent$:$
\begin{itemize}

\item[(1)] $J$ is multiplicatively closed.

\item[(2)] There is $\NN \succ \MM$ such that $\gcis(\MM,\NN) = J$.

\item[(3)] $\MM$ has a minimal extension $\NN$ such that $\gcis(\MM,\NN) = J$ and $\NN$ does not fill $J$.

\end{itemize}}

\bigskip

{\it Proof}. The equivalence $(1) \Longleftrightarrow (2)$ is from Paris-Mills \cite{pm}. It's trivial that $(3) \Longrightarrow (2)$. Even though $(1) \Longrightarrow (3)$ follows from \cite[Th.~7.6]{bases}, we indicate how to prove $(1) \Longrightarrow (3)$ since the ideas come up again. Given that $J$ is a multiplicatively closed cut of $\MM$, how do we get $\NN \succ \MM$ satisfying the following?
\smallskip

\begin{itemize}

\item[]
\begin{itemize}
\item[(M1)] $\NN$ is a minimal extension of $\MM$.
\item[(M2)]   $\gcis(\MM,\NN) = J$.

\item[(M3)]  $\NN$ does not fill $J$.
\end{itemize}
\end{itemize}
\smallskip
 Such an $\NN$ can be constructed from 
 a decreasing sequence $X_0 \supseteq X_1 \supseteq X_2 \supseteq \cdots$ of infinite sets in $\Def(\MM)$ such that:
 \smallskip

\begin{itemize}

\item[]
\begin{itemize}

\item[(S1)] For every $\MM$-definable function $f : M \into M$, there is $i < \omega$ such that $f$ is either constant or one-to-one on $X_i$.

\item[(S2)] For  every $d \in M \backslash J$ there is $i < \omega$ such that $\MM \models |X_i| \leq d$.

\item[(S3)] For every  $\MM$-definable one-to-one function $f : M \into M$, there are $i < \omega$ and  $d \in M \backslash J$ such that $f(x) \geq d$ for every $x \in X_i$.

\end{itemize}

\end{itemize}
\smallskip

This sequence determines a unique complete, nonprincipal $1$-type $p(x)$ over $\MM$; to wit: if $\theta(x)$ is an ${\mathcal L}(M)$-formula that defines $A$ in $\MM$, then $\theta(x) \in p(x)$ iff there is $i < \omega$ such that $X_i \subseteq A$. 
To see that $p(x)$ is  nonprincipal and complete, we just need  each $X_i$ to be infinite and the following, which is a consequence of  (S1): 

\smallskip

\begin{itemize}

\item[]
\begin{itemize}
\item[(S0)] For every $A \in \Def(M)$, there is $i < \omega$ such that $X_i \subseteq A$ or $X_i \cap A = \varnothing$.
\end{itemize}
\end{itemize}
\smallskip
Having $p(x)$, we let $\NN \succ \MM$ be generated by  $c$ over $\MM$, where $c$ realizes $p(x)$. Alternatively, $\NN$ is generated over $\MM$ by an element $c$, where $c \in X_i^\NN$ for all $i < \omega$. Incidentally, $\{c\} = \bigcap_{j< \omega}X_j^\NN$ by Ehrenfeucht's Lemma \cite{eh} (or \cite[Th.~1.7.2]{ksbook}).

Each of  (M1) -- (M3) is satisfied. Only (S1) is needed to guarantee (M1). The $\gcis(\MM,\NN) \subseteq J$ half of (M2) follows from (S2).  The other half of (M2), $\gcis(\MM,\NN) \supseteq J$, follows from (S3), or even from the following consequence of  (S3): 
\smallskip

\begin{itemize}
\item[]
\begin{itemize}

 \item[(S3$'$)] For every $i < \omega$, there is $d \in M \backslash J$ such that $\MM \models |X_i| \geq d$. 
 \end{itemize}
 \end{itemize}
 \smallskip
 Note that (S3) is sufficiently stronger than (S3$'$) to imply that (M3) is satisfied. 

We next say a little more about the sequence of $X_i$'s. Each $X_i$ is large, where we say that $X$ is {\em large} if $X$ is $\MM$-finite and $\MM \models |X| = c$ for some $c \in M \backslash J$. Corresponding to (S1) -- (S3) respectively, there are the following facts concerning a large $X$:
\smallskip

\begin{itemize}
\item[]
\begin{itemize}

 \item[(F1)]  If   $f : X \into M$ is $\MM$-definable, then there is   large $Y \subseteq X$ such that $f$ is either constant or one-to-one on $Y$.

\item[(F2)] If  $d \in M \backslash J$, then there is  large $Y \subseteq X$ such that $\MM \models |Y| \leq d$.

\item[(F3)] If  $f : X \into M$ is one-to-one and $\MM$-definable, then there are large $Y \subseteq X$ and  $d \in M \backslash J$ such that $\MM \models \forall x \in Y[f(x) \geq d]$. 
\end{itemize}

\end{itemize}
\smallskip
Easily, (F2) is trivial, (F1) follows from {\sf PHP2}, and (F3) from {\sf PHP1}.  The countability of $\MM$ allows us to get (S1) -- (S3) from (F1) -- (F3). \qed

\bigskip

{\sc Remark}: In the previous proof, (F3) is a consequence of (F1). To see this, let $d$ be the least such that $|\{x \in X : f(x) > d\}| \leq |\{x \in X : f(x) \leq d\}|$. Then let $Y = \{x \in X : f(x) \geq d\}$.  
 \bigskip

{\sc Theorem 1.2}: (Probably  known) {\em Suppose that $J$ is a cut of countable $\MM$. The following are equivalent$:$
\begin{itemize}

\item[(1)] $J$ is exponentially closed.

\item[(2)] There is $\NN \succ \MM$ such that $\gcis(\MM,\NN) = J$ and $\cod(\NN / J) = \cod(\MM / J)$.

\item[(3)] $\MM$ has a minimal extension $\NN$ such that $\gcis(\MM,\NN) = J$, $\cod(\NN / J) = \cod(\MM / J)$
 and $\NN$ does not fill $J$.
 
 \end{itemize}} 
 
 \bigskip
 
 {\it Proof}. It suffices to prove only $(2) \Longrightarrow (1)$ and $(1) \Longrightarrow (3)$, since $(3) \Longrightarrow (2)$ is trivial.
 
  For $(2) \Longrightarrow (1)$, assume that $J$ is not exponentially closed, but that  $\NN \succ \MM$ and $\gcis(\MM,\NN) = J$. Let $a \in J$ be such that $2^a \not\in J$. Let $b \in N \backslash M$ be such that $b < 2^a$.  
 
For $x \in M$, let $\Bin(x) \in \Def(\MM)$ be uniquely defined by  $\MM \models x = \sum\{2^e : e \in \Bin(x)\}$. Thus,  $\cod(\MM / J) = \{J \cap \Bin(x) : x \in M\}$.
Then $\Bin(b)  \subseteq [0,a) \subseteq J$ and  $\Bin(b)  \not\in \Def(\MM)$. 
 (Note that this implies not only that $\cod(\NN / J) \neq \cod(\MM / J)$, but something apparently stronger:  
 $\Def(\NN) \cap {\mathcal P}(J) \neq \Def(\MM) \cap {\mathcal P}(J)$.)

 For $(1) \Longrightarrow (3)$, we follow the proof of $(1) \Longrightarrow (3)$ of Theorem~1.1, but now we  need, in addition to (M1) -- (M3), the following: 
 \smallskip 
 
 \begin{itemize}
 \item[]
 \begin{itemize}
 
 \item[(M4)] $\cod(\NN / J) = \cod(\MM / J)$.

 \end{itemize}
 \end{itemize}
 \smallskip
 This will follow if the sequence of $X_i$'s  satisfies, in addition to (S1) -- (S3), the following:
 \smallskip
 \begin{itemize}
 \item[]
 \begin{itemize}
 
 \item[(S4)] For every  $\MM$-definable function $g : M \into \Def(\MM)$, there are $i < \omega$ and $d \in M \backslash J$ such that for all $x,y \in X_i$, $g(x) \cap [0,d) = g(y) \cap [0,d)$. 
 \end{itemize}
 \end{itemize}
 \smallskip
 Since $J$ is exponentially closed, we can get the sequence to satisfy (S4). This is a consequence of the following fact concerning a large $X$:
 \smallskip
 
 \begin{itemize}
 \item[]
 \begin{itemize}
 
 \item[(F4)] There is $b \in M \backslash J$ such that if $g : X \into {\mathcal P}([0,b))$ is $\MM$-definable, then then there is large $Y \subseteq X$ such that $g$ is constant on $Y$.
 \end{itemize}
 \end{itemize} 
 \smallskip
 But (F4) is a consequence of (F1). To see why, we first consider the following strengthening of (F1):
 \smallskip
 \begin{itemize}
 \item[]
 \begin{itemize}
 
 \item[(F$1'$)] There is $d \in M \backslash J$ such that if $f : X \into M$ is $\MM$-definable, then there is $Y \subseteq X$ such that $\MM \models |Y| = d$ and $f$ is either constant or one-to-one on $Y$.
 \end{itemize}
 \end{itemize} 
 \smallskip
 (F$1'$) is a consequence of (F1). Let $e \in M$ be big enough so that there is $\{f_i : i < e\}$ consisting of all $\MM$-definable $f : X \into X$. For each $i < e$, let $d_i$ be the largest for which there is $Y \subseteq X$, $|Y| = d_i$ and $f_i$ is either constant or one-to-one on $Y$. By (F1), each $d_i \not\in J$. Since $J \not\in \Def(\MM)$, there is $d \in M \backslash J$ such that $d \leq d_i$ for all $i < e$.
 To see that this $d$ works, consider some $\MM$-definable $f : X \into M$ and then let $i < e$ be such that for $x,y \in X$, $f_i(x) = f_i(y)$ iff $f(x) = f(y)$. Let $Y \subseteq X$ be such that $|Y| = d$ and  $f_i$ is constant or one-to-one on $Y$. Then,   $f$ is constant or one-to-one on $Y$.
 
 Now to get (F4), let $d$ be as in (F$1'$) and let $b \in M \backslash J$ be such that $2^b < d$. 
 Letting $Y \subseteq X$ be such that $|Y|=d$ and $g$ is either constant or one-to-one on $Y$, we see that $g$ is, in fact, constant on $Y$. \qed
 
 \bigskip
 
 The following proposition, whose straightforward proof is omitted, asserts that each of the properties (M1), (M3),  (M4) is characterized by an ${\mathcal L}^*$-sentence.
 
  \bigskip

{\sc Proposition 1.3}: {\em There are ${\mathcal L}^*$\hspace{-1pt}-sentences~$\sigma_1, \sigma_3,\sigma_4$ such that whenever $J$ is a cut of $\MM$, $\NN \succ \MM$ and $\gcis(\MM,\NN) = J$, then$:$

$(1)$ $(\NN,M) \models \sigma_1 \Longleftrightarrow \NN$ is a minimal extension of $\MM$.

$(3)$ $(\NN,M) \models \sigma_3 \Longleftrightarrow \NN {\mbox{ does not fill }} J.$

$(4)$ $(\NN,M) \models \sigma_4 \Longleftrightarrow \cod(\NN / J) = \cod(\MM / J)$. \qed}

\bigskip

%
%

\subsection{Lattices} We say a few words about the lattices $\Lt(\NN / \MM)$, which are defined in the introduction. Let
$$
\Lt_0(\NN / \MM) = \{\KK \in \Lt(\NN / \MM) : \KK {\mbox{ is finitely generated over }}\MM\}.
$$
 It is possible that 
$\Lt_0(\NN / \MM)$ is not a sublattice of $\Lt(\NN / \MM)$ since there may be $\KK_1,\KK_2 \in \Lt_0(\NN / \MM)$ with  $\KK_1 \wedge \KK_2 \not\in \Lt_0(\NN / \MM)$. However, $\Lt_0(\NN / \MM)$ is a $\vee$-subsemilattice of $\Lt(\NN / \MM)$, and we typically think of $\Lt_0(\NN / \MM)$ as a $\vee$-semilattice.  If $\Lt(\NN / \MM)$ is finite, then  $\Lt_0(\NN / \MM)$ $= \Lt(\NN / \MM)$.  It is observed in Lemma~7.2 of \cite{bases}  that if $\MM \prec_{\sf cf} \NN$, then 
$\Lt_0(\NN / \MM)$ is interpretable in $(\NN,M)$. Moreover, as one sees from the proof,  there is  uniform interpretation whenever $\MM \prec_{\sf cf} \NN$. 

It is asked (\cite[Question~7.3]{bases}) if $\Lt(\NN_1 / \MM_1) \equiv \Lt(\NN_2 / \MM_2)$ whenever 
$\MM_1 \prec_{\sf cf} \NN_1$, $\MM_2 \prec_{\sf cf} \NN_2$ and $(\NN_1,M_1) \equiv (\NN_2, M_2)$. We show that this question has negative answer by giving a 
counter-example. Consider any nonstandard countable $\MM_1$. Then there are $\NN_1 \succ_{\sf cf} \MM_1$ such that $\Lt(\NN_1/\MM_1)$ is a chain of length $\omega+1$. Then $\Lt_0(\NN_1/\MM_1)$ is a chain of length $\omega$. Now consider a countable, recursively saturated $(\NN_2,M_2) \equiv (\NN_1,M_1)$. 
Then, $\Lt_0(\NN_2 / \MM_2)$ is a chain whose order-type $\omega + \ZZ \cdot \eta$, where $\ZZ$ is the order-type of the integers and $\eta$ is the order-type of the rationals. It is clear that $\Lt(\NN_1 / \MM_1)$ has only one limit point but $\Lt(\NN_2 / \MM_2)$ has more than one (in fact, uncountably many). Thus, 
$\Lt(\NN_1 / \MM_1) \not\equiv \Lt(\NN_2 / \MM_2)$.  

For a set $A$, we let $\eq(A)$ be the set of equivalence relations on $A$. If $\Theta \in \eq(A)$, 
then a $\Theta$-class is an equivalence class of $\Theta$. Furthermore, $\Theta$ is {\em trivial} if $A$ is the only $\Theta$-class, and $\Theta$ is {\em discrete} if each $\Theta$-class is a singleton. We sometimes think of $\eq(A)$ as a lattice, whose minimal element  is  the discrete equivalence relation $\bbzero_A$ and whose maximal element is  the trivial equivalence relation $\bbone_A$. 

If $L$ is a finite lattice, then a {\it representation} of $L$ is a one-to-one function $\alpha : L \into \eq(A)$ such that  $\alpha(0_L) = \bbone_A$, $\alpha(1_L) = \bbzero_A$ and $\alpha(r \vee s) = \alpha(r) \cap \alpha(s)$. If $\alpha : L \into \eq(A)$ is a representation and $B \subseteq A$, then 
$\alpha |B : L \into \eq(B)$ where $(\alpha | B)(r) = \alpha(r) \cap B^2$  for every $r \in L$. The representation $\alpha$ is 0-CPP if, for each $r \in L$,  $\alpha(r)$ does not have exactly 2 equivalence classes. Recursively, $\alpha$ is $(n+1)$-CPP if whenever $\Theta \in \eq(A)$, then there is $B \subseteq A$ such that $\alpha |B$ is an $n$-CPP representation of $L$ and for some $r \in L$, $\alpha(r) \cap B^2 = \Theta \cap B^2$. These definitions can be found in \cite[Chap.~4.5]{ksbook} and their connection with 
interstructure lattices can also be found there. We assume that the reader has some familiarity with the material in Chapter~4 of \cite{ksbook}. 

Suppose that $K$ and $L$ are finite lattices. Thus, both $1_K$ and $0_L$ exist. We define their {\em linear sum} $K \oplus L$ to lattice that is the disjoint union of $K$ and $L$ (except that we identify $1_K$ and $0_L$) such that both $K$ and $L$ are sublattices of $K \oplus L$ and if $r \in K$ and $s \in L$, then $r \leq s$ in $K \oplus L$. For example, ${\mathbf 2} \oplus {\mathbf 2} \cong {\mathbf 3}$.

 \bigskip

%
%

\section{Humble Extensions} Before introducing humble extensions, we begin this section with a simple observation concerning finitely generated extensions. 

 \bigskip
 
 {\sc Proposition 2.1}: {\em Suppose that $\MM \prec \NN$ and that $a$ generates $\NN$ over $\MM$. Let $\MM \prec \NN_0 \preccurlyeq \NN$, where $\NN_0$ is finitely generated over $\MM$. Then  there is $b \leq a$ that generates $\NN_0$ over $\MM$.}
 
 \bigskip
 
 {\it Proof}. Let $a_0$ generate $\NN_0$ over $\MM$. Then there is an $\MM$-definable function $t_0 : M \into M$ such that $\NN \models a_0 = t_0(a)$. Using recursion in $\MM$, define the function $t : M \into M$ so that 
  \begin{equation*}
 t(x) =
 \begin{cases}
 x,  \mbox{ if } t_0(z) \neq t_0(x) \mbox{ for all } z< x, \\
 
 t(z), {\mbox{ where }} z = \min\{y < x : t_0(y) = t_0(x)\}, \mbox{ otherwise.}
 \end{cases}
 \end{equation*}
 Clearly, $t$ is $\MM$-definable and $\MM \models \forall x[t(x) \leq x]$. Let $b = t^\NN(a)$, so that $\NN \models b \leq a$.  We will show that $b$ generates $\NN_0$ over $\MM$.   This can be accomplished by showing that there are $\MM$-definable functions $t_1, t_2$ such that 
 $$
 \NN \models t_1(a_0) = b \wedge t_2(b) = a_0.
 $$
 Note that $t_1, t_2$ need not be total, but it should be that $a_0 \in \dom(t^\NN_1)$ and $b \in \dom(t^\NN_2)$.
 
  Let $R, R_0$ be the ranges of $t, t_0$, respectively. Let $t_1 : R \into M$ and $ t_2 : R_0 \into M$ be such that   
 \begin{equation} \tag{$*$}
 t_1(t(x)) = t_0(x) \mbox{ and }t_2(t_0(x)) = t(x)
 \end{equation} 
 for each $x \in M$. 
Observe that both $t_1$ and $t_2$ are well defined since for any $x,y \in M$, $t(x) = t(y)$ iff $t_0(x) = t_0(y)$. Thus, both $t_1,t_2$ are $\MM$-definable. It follows from $(*)$ that
$$
t_1^\NN(b) = t_1^\NN(t^\NN(a)) = t^\NN_0(a) = a_0
$$
and
$$
t^\NN_2(a_0) =  t^\NN_2(t^\NN_0(a)) = t^\NN(a) = b,
$$
thereby proving that $b$ generates $\NN_0$ over  $\MM$.  \qed
 
 \bigskip
 
  Suppose that $\MM$ is nonstandard. Thus, $\omega$ is the {\em standard} cut of $\MM$. In this section, we will be concerned with elementary extensions $\NN \succ \MM$ that are generated over $\MM$ by an element $c$ that fills the standard cut of $\MM$, meaning that $\omega = \{a \in M : \NN \models a < c\}$. We refer to such  extensions as a  {\bf humble} extensions. If $\NN$ is a humble extension of $\MM$, then $\MM \prec_{\sf cf} \NN$ and $|N| = |M|$.
  
  \bigskip
  
  {\sc Corollary 2.2}: {\em Suppose that $\NN$ is a humble extension of $\MM$, and $\MM \prec \NN_0 \preccurlyeq \NN$, where $\NN_0$ is finitely generated over $\MM$. Then $\NN_0$ is a humble extension of $\MM$.} 
  
  \bigskip
  
  {\it Proof}. Let $a \in N$ fill the standard cut of $\MM$ and generate $\NN$ over $\MM$. Proposition~2.1 implies that there is $b \leq a$ that generates $\NN_0$ over $\MM$. Obviously, $b$ also fills the standard cut of $\MM$. \qed 
  
  \bigskip

 Recall that the standard system of nonstandard $\MM$ is $\ssy(\MM) = \cod(\MM / \omega)$, which is a Boolean subalgebra of ${\mathcal P}(\omega)$ and  which contains all  finite and cofinite subsets of $\omega$. A filter ${\mathcal F}$ of $\ssy(\MM)$ is  {\em free} (or {\em Fr\'echet}) if it contains all cofinite subsets of $\omega$.  An ultrafilter ${\mathcal U}$ of $\ssy(\MM)$ is nonprincipal iff it is free. 
   
Humble extensions of a nonstandard $\MM$ can be obtained in the following way. Let ${\mathcal U}$ be a nonprincpal ultrafilter of $\ssy(\MM)$. Then ${\mathcal U}$ determines a complete 1-type  $p(x)$ 
over $\MM$ as follows: if $\theta(x)$ is an ${\mathcal L}(\MM)$-formula, then $\theta(x) \in p(x)$ iff there is $A \in {\mathcal U}$ such that $\MM \models \theta(a)$ for all $a \in A$. Let $\NN$ be generated over $\MM$ by an element $c$ that realizes $p(x)$. This $\NN$ is a humble extension of $\MM$ generated by $c$, which fills the standard cut. 

Alternatively,  the extension $\NN$  can be viewed as a kind of ultrapower, the well known construction of which we now describe.   Let $T$ be the set of  $\MM$-definable functions $t : M \into M$. Define the equivalence relation $\sim$ on $T$ so that if $t_1,t_2 \in T$, then $t_1 \sim t_2$ iff $\{n < \omega : t_1(n) = t_2(n)\} \in {\mathcal U}$. Let $[t]$ be the equivalence class to which $t$ belongs. 
By a mixture of theorems of Skolem \cite{ts} and Lo\'{s} \cite{los}, there is $\NN$ such that $N = \{[t] : t \in T\}$ and for any ${\mathcal L}$-formula 
$\varphi(x_0, \dots, x_{j-1})$ and $t_0, \ldots, t_{j-1} \in T$, 
$$
\NN \models \varphi([t_0], \ldots, [t_{j-1}]) {\mbox{ iff }} \{n < \omega : \MM \models \varphi(t_0(n), \ldots, t_{j-1}(n))\} \in {\mathcal U}.
$$
By identifying $a \in M$ with $[t_a]$, where $t_a : M \into M$ is the $\MM$-definable function that is constantly $a$, we get that $\MM \prec \NN$.   Thus, we use the suggestive notation 
$$
 \NN = \MM^\omega / {\mathcal U}
$$
when $\NN$ is constructed in either of these ways. Notice that $[id_M]$ generates $\NN$ over $\MM$, where $id_M : M \into M$ is the identity function.

Conversely, suppose that  $\NN$ is a humble extension of $\MM$ generated over $\MM$ by $c$, which fills the standard cut. Let 
$$
 {\mathcal U}_c = \{X \cap \omega \in \ssy(M) :  X \in \Def(\MM) \mbox{ and }  a \in X^\NN\}.
 $$
 Then, ${\mathcal U}_c$ is a nonprincipal ultrafilter of $\ssy(\MM)$ and $ \MM^\omega / {\mathcal U}_c \cong \NN$. 
 In fact, there is an isomorphism $f $ such that if $t : M \into M$ is $\MM$-definable, then $f([t]) = t^\NN(c)$, so that, in particular, $f([id_M]) = c$.
 
 Summarizing:
 
 \bigskip
 
 {\sc Proposition 2.3}: {\em Suppose that $\MM$ is nonstandard.} 
 
  $(1)$ {\em If ${\mathcal U}$ is a nonprincipal ultrafilter of $\ssy(\MM)$, then  $\NN = \MM^\omega/{\mathcal U}$
 is a humble extension of $\MM$ generated by $[id_M]$.} 
 
 $(2)$ {\em If $\NN$ is a humble extension of $\MM$ generated by $c$ over $\MM$, then there is an isomorphism $f : \NN \into \MM^\omega/{\mathcal U}_c$ that fixes $M$  such that $f(c) = [id_M]$.}
 \qed
 
 \bigskip

 \hyphenation{ultra-filters}
 
  The well known Rudin-Keisler quasi-order $\leq_{\sf RK}$ on ultrafilters over $\omega$ has a natural generalization to those ultrafilters of $\ssy(\MM)$.  Recall its definition: If ${\mathcal U},{\mathcal V}$ are ultrafilters over $\omega$, then ${\mathcal U} \leq_{\sf RK} {\mathcal V}$ iff there is $f : \omega \into \omega$ such that whenever $A \in {\mathcal U}$, then $f^{-1}(A) \in {\mathcal V}$. We generalize $\leq_{\sf RK}$ from ${\mathcal P}(\omega)$ to $\ssy(\MM)$.
  
  If  ${\mathcal U}$, ${\mathcal V}$ are ultrafilters of $\ssy(\MM)$, then we define ${\mathcal U} \leq_{\sf RK}^\MM {\mathcal V}$ if there is an $\MM$-definable $t : M \into M$  such that 
  whenever $X \in \Def(\MM)$ and $X \cap \omega \in {\mathcal U}$, then  $t^{-1}(X) \cap \omega \in {\mathcal V}$. For such a $t$, we will say that $t$ {\em demonstrates} that ${\mathcal U} \leq_{\sf RK}^\MM {\mathcal V}$. One easily checks that  $\leq_{\sf RK}^\MM$ is a quasi-order $-$ in other words, if ${\mathcal U}, {\mathcal V}, {\mathcal W}$ are ultrafilters of $\ssy(\MM)$ and 
  ${\mathcal U} \leq^\MM_{\sf RK}   {\mathcal V} \leq^\MM_{\sf RK}   {\mathcal W}$, then ${\mathcal U} \leq^\MM_{\sf RK}   {\mathcal W}$. Moreover, if $s$ demonstrates that ${\mathcal U} \leq^\MM_{\sf RK}   {\mathcal V}$ and $t$ demonstrates that   ${\mathcal V} \leq^\MM_{\sf RK}   {\mathcal W}$, then 
  $s \circ t$ demonstrates that ${\mathcal U} \leq^\MM_{\sf RK}   {\mathcal W}$.
  
  By the next proposition, there is a simpler definition of $\leq^\MM_{\sf RK}$ if $({\mathbb N}, \ssy(\MM) ) \models {\sf ACA_0}$.

  \bigskip
  
  {\sc Proposition 2.4}: {\em If $\MM$ is nonstandard, $({\mathbb N}, \ssy(\MM) ) \models {\sf ACA_0}$ and 
 ${\mathcal U}$, ${\mathcal V}$ are ultrafilters of $\ssy(\MM)$, then following are equivalent.
  
  \begin{itemize}
  
  \item[$(1)$] ${\mathcal U} \leq_{\sf RK}^\MM {\mathcal V};$ 
  
 \item[$(2)$] 
  there is $f : \omega \into \omega$, coded in $\ssy(\MM)$, such that  $f^{-1}(A) \in {\mathcal V}$ whenever $A \in {\mathcal U}$. 
  
  \end{itemize}}
  
  \bigskip 
  
  {\sc Remark}: Note that $(\mathbb{N}, \ssy(\MM)) \models {\sf ACA}_0$ iff $\omega$ is a {\em strong} cut of~$\MM$.  
  
  \bigskip
  
  {\it Proof}. Suppose that  ${\mathcal U}, {\mathcal V}$ are ultrafilters of $\ssy(\MM)$, where $\MM$ is nonstandard and $({\mathbb N}, \ssy(\MM) ) \models {\sf ACA_0}$. 
. 
  
  \smallskip
  
  $(1) \Longrightarrow (2)$: Suppose that ${\mathcal U} \leq_{\sf RK}^\MM {\mathcal V}$. Let $t : M \into M$ demonstrate that ${\mathcal U} \leq_{\sf RK}^\MM {\mathcal V}$. Let $f = t \cap \omega^2$, so that $f : D \into \omega$ for some $D \subseteq \omega$.
  Since $({\mathbb N}, \ssy(\MM) ) \models {\sf ACA_0}$, then $D \in \ssy(\MM)$ and $R= f[D] \in \ssy(\MM)$. We claim that $R \in {\mathcal U}$. 
  
  If $R \not\in {\mathcal U}$, then  $\omega \backslash R \in {\mathcal U}$. Let $X \in \Def(\MM)$ be such that $X \cap \omega = \omega \backslash R$. Then  $\varnothing = t^{-1}(X) \cap \omega \in {\mathcal V}$, which is a contradiction showing that $R \in {\mathcal U}$.

  In order to prove (2), suppose that $A \in {\mathcal  U}$, intending to show that $f^{-1}(A) \in {\mathcal V}$. Let $X \in \Def(\MM)$ be such that $A = X \cap \omega$. It suffices to show that $f^{-1}(A\cap R) \in {\mathcal V}$. For a contradiction, suppose that $f^{-1}(A \cap R) \not\in {\mathcal V}$. Since $D = f^{-1}(R) \in {\mathcal V}$, then $f^{-1}(R \backslash A) \in {\mathcal V}$.
  But $t^{-1}(M \backslash X) \cap \omega = 
f^{-1}(\omega \backslash A) \in {\mathcal V} \supseteq f^{-1}(R \backslash A)$, so $t^{-1}(M \backslash X) \cap \omega \in {\mathcal V}$. This implies $t^{-1}(X) \cap \omega \not\in {\mathcal V}$, which, in turn, implies that $X \cap \omega \not\in {\mathcal U}$, a contradiction. 

Of course, this $f$ may not have $\omega$ as its domain. If that is the case, replace it with an extension of itself that is coded in $\ssy(\MM)$ having domain of $\omega$.

\smallskip

$(2) \Longrightarrow (1)$: (This implication does not require that $({\mathbb N},\ssy(\MM)) \models {\sf ACA}_0$.) Let $f$ be as in $(2)$. Since $f$ is coded in $\ssy(\MM)$, there is an $\MM$-definable $T \subseteq M^2$ such that $T \cap \omega^2 = f$. Let 
$$
t = \{\langle x,y \rangle \in M^2 : y = \min(\{0\} \cup \{z \in M : \langle x,z \rangle \in T\})\}.
$$
 Then $f = t \cap \omega^2$ and $t : M \into M$. To see that $t$ demonstrates that 
 ${\mathcal U} \leq_{\sf RK}^\MM {\mathcal V}$, 
 suppose that $X \in \Def(\MM)$ is such that $X \cap \omega \in {\mathcal U}$. Then,   $t^{-1}(X) \cap \omega = f^{-1}(X \cap \omega) \in {\mathcal V} $.  \qed
 
 \bigskip
 
  As said earlier, $\leq_{\sf RK}^\MM$ is a generalization of $\leq_{\sf RK}$. The next corollary makes this precise. 
 
 \bigskip
 
 {\sc Corollary 2.5}: {\em Suppose that $\ssy(\MM) = {\mathcal P}(\omega)$ and that ${\mathcal U},{\mathcal V}$ are ultrafilters over $\omega$. Then, ${\mathcal U} \leq_{\sf RK} {\mathcal V}$ iff 
 ${\mathcal U} \leq^\MM_{\sf RK} {\mathcal V}$.}
 
 \bigskip
 
 {\it Proof}. If $\ssy(\MM) = {\mathcal P}(\omega)$, then $({\mathbb N},\ssy(\MM)) \models {\sf ACA}_0$. \qed
 
  \bigskip
 
  {\sc Lemma 2.6}: {\em Let $\MM$ be  nonstandard, and let ${\mathcal U}, {\mathcal V}$ be  ultrafilters of $\ssy(\MM)$.  Then, 
  ${\mathcal U} \leq^\MM_{\sf RK} {\mathcal V}$ iff 
 there is an elementary embedding of $\MM^\omega/{\mathcal U}$ into $\MM^\omega/{\mathcal V}$ fixing~$M$.}
 
  \bigskip
 
  {\it Proof}. Let $\MM$,  ${\mathcal U}$ and ${\mathcal V}$ be as given. Let $\NN_1 = \MM^\omega/{\mathcal U}$ and $\NN_2 = \MM^\omega/{\mathcal V}$.
  
  \smallskip
  
  $(\Longrightarrow)$: Letting $T$ be the set of $\MM$-definable functions $t : M \into M$, we see that ${\mathcal U}$ and ${\mathcal V}$ determine, respectively,  the two equivalence relations $\sim_1$ and $\sim_2$ in $\eq(T)$. 
 For $t \in T$, we let $[t]_1, [t]_2$ be  the respective  equivalence classes to which $t$ belongs.
 
  Suppose that $F : M \into M$ demonstrates that ${\mathcal U} \leq^\MM_{\sf RK} {\mathcal V}$. Thus,  $F$ is  $\MM$-definable, 
  $N_1 = \{[t]_1 : t \in T\}$ and  $N_2 = \{[t]_2 : t \in T\}$. Let $g : N_1 \into \NN_2$ be such that $g([t]_1) = [t \circ F]_2$. 
 
 First, we show that $g$ is well-defined: if $t \sim_1 s$, then $t \circ F \sim_2 s \circ F$.
 \begin{align*}
  t \sim_1 s 
                  & \Longrightarrow \{n \in \omega :  t(n) = s(n)\}  \in {\mathcal U} \\
                  & \Longrightarrow F^{-1}(\{x \in M : \MM \models t(x) = s(x)\}) \cap \omega \in {\mathcal V} \\
                  & \Longrightarrow \{m < \omega : \MM \models t(F(m)) = s(F(m))\} \in {\mathcal V} \\
                  & \Longrightarrow t \circ F \sim_2 s \circ F .
     \end{align*}
     
Next, we show that $g$ is an elementary embedding. 
 \begin{align*}
   \NN_1 \models \varphi\big([t_0]_1, \ldots, [t_{j-1}]_1\big)
                                           & \Longrightarrow   \{n < \omega : \MM \models \varphi\big(t_0(n), \ldots,t_{j-1}(n)\big)\} \in {\mathcal U} 
                                           \end{align*}
                                           \vspace{-24pt}
                                           \begin{align*}
                                            & \Longrightarrow  F^{-1}( \{x \in M   : \MM \models \varphi\big(t_0(x), \ldots, t_{j-1}(x)\big)\}) \cap \omega \in {\mathcal V} \\ 
                                           & \Longrightarrow  \{m < \omega : \MM \models \varphi\big(t_0(F(m)), \ldots, t_{j-1}(F(m))\big)\} \in {\mathcal V} \\ 
                                            & \Longrightarrow \NN_2 \models \varphi\big([t_0 \circ F]_2, \ldots, [t_{j-1} \circ F]_2\big) \\
                                            & \Longrightarrow \NN_2 \models \varphi\big(g([t_0]_1), \ldots, g([t_{j-1}]_1)\big). 
 \end{align*} 
It is clear that $g$ fixes $M$. 

\smallskip
 
$(\Longleftarrow)$: Let $g$ be an elementary embedding of $\NN_1$ into $\NN_2$ fixing $M$.  Let $a_1 = [id_M]_1 \in N_1$ and $a_2 = [id_M]_2 \in N_2$.   Obviously, both $a_1$ and $a_2$ fill the standard cut of $\MM$.   Let $b = g(a_1) \in N_2$. Therefore, 
${\mathcal V} = \{X \cap \omega : X \in \Def(\MM)$ and $b \in X^{\NN_2}\}$, so we can safely assume that $g$ is the identity function. Hence, $b = a_1$.

Let $t : M \into M$ be an $\MM$-definable function such that $\NN_2 \models a_1 = t(a_2)$.    We claim that $t$ demonstrates that ${\mathcal U} \leq_{\sf RK}^\MM {\mathcal V}$. 

To prove the claim, consider some 
$X \in \Def(\MM)$ and $X \cap \omega \in {\mathcal V}$. Thus, $a_1 \in X^{\NN_2}$. 
Let $Y = t^{-1}(X)$. Then, $Y \in \Def(\MM)$ and $a_2 \in t^{-1}(\{a_1\}) \subseteq Y$. Hence, $Y \cap \omega \in {\mathcal V}$, proving the claim. \qed

\bigskip

If ${\mathcal U}, {\mathcal V}$ are ultrafilters of $\ssy(\MM)$, then ${\mathcal U}\cong^\MM_{\sf RK} {\mathcal V}$ if  ${\mathcal U}\leq^\MM_{\sf RK} {\mathcal V} \leq^\MM_{\sf RK} {\mathcal U}$.

\bigskip

{\sc Corollary 2.7}: {\em Let $\MM$ be  nonstandard, and let ${\mathcal U}, {\mathcal V}$ be  ultrafilters of $\ssy(\MM)$.  Then, 
  ${\mathcal U} \cong^\MM_{\sf RK} {\mathcal V}$ iff 
 there is an isomorphism of $\MM^\omega/{\mathcal U}$ and $\MM^\omega/{\mathcal V}$ fixing~$M$.}
 
 \bigskip
 
 {\it Proof}. Let $\NN_1 = \MM^\omega/{\mathcal U}$ and $\NN_1 = \MM^\omega/{\mathcal V}$.
Since both $\NN_1, \NN_2$ are finitely generated over $\MM$ and there there are elementary embeddings $f_1 : \NN_1 \into \NN_2$ and $f_2 : \NN_2 \into \NN_1$ fixing $\MM$, the it follows from Ehrenfeucht's Lemma \cite{eh} (or \cite[Th.~1.7.2]{ksbook}), that there  is an isomorphism of $\NN_1$ and $\NN_2$ fixing~$M$. \qed

  \bigskip
  
  {\sc Lemma 2.8}: {\em Suppose that $\MM$ is nonstandard and countable. There is $\NN \succ \MM$ that is a minimal, humble extension.}
  
  \bigskip
  
  {\it Proof}. We construct $\NN$ from a decreasing sequence $X_0 \supseteq X_1 \supseteq X_2 \supseteq \cdots$ of sets in $\Def(\MM)$ such that:
  
  \begin{itemize}  
  
  \item[(1)] For each $i < \omega$, $X_i \cap \omega$ is infinite.
  
  \item[(2)] For each $\MM$-definable $f : M \into M$, there is $i < \omega$ such that $f$ is either constant or one-to-one on $X_i$.
  
  \end{itemize}
  Having this sequence, we let ${\mathcal U} = \{A \in \ssy(\MM): A \subseteq X_i$ for some $i < \omega\}$, and then let $\NN = \MM^\omega / {\mathcal U}$. One easily checks that, first, ${\mathcal U}$ is an ultrafilter of $\ssy(\MM)$ and, second, that $\NN$ is a minimal, humble extension of $\MM$.
  
  To get the sequence, we make use of the countability of $\MM$ and the following fact.
  
  \smallskip
  
  \begin{quote} 
   Let $f : M \into M$ be
  $\MM$-definable and $X \in \Def(\MM)$ be such that $X \cap \omega$ is infinite. Then there is $Y \subseteq X$ such that $Y \in \Def(\MM)$, $Y \cap \omega$ is infinite and $f$ is either constant or one-to-one on $Y$.
  
  \end{quote}
  To prove this fact, let $f$ and $X$ be as indicated. If there is $a \in M$ such that $\{n \in X : f(n) = a\}$ is infinite, then let $Y = \{x \in X : f(x) = a\}$. If there is no such $a$, then let 
  $Y = \{x \in X : f(x) \neq f(y)$ for all $y \in X \cap [0,x)\}$. \qed

  \bigskip
  
   We make a short digression about humble extensions and their interstructure lattices.

 \bigskip
 
 {\sc Lemma 2.9}: {\em For countable nonstandard $\MM$ and finite lattice $L$, the following are equivalent$:$
 
 \begin{itemize} 
 
 \item [(1)] $L$ has $n$-CPP representations for each $n < \omega$.
 
 \item[(2)] There is a humble extension $\NN \succ \MM$ such that \\ $\Lt(\NN / \MM) \cong {\mathbf 2} \oplus L$.
 
 \end{itemize}}
 
 \bigskip
 
 {\it Proof}. Let $\MM$ be countable and nonstandard, and let $L$ be a finite lattice.
  
 \smallskip
 
 $(2) \Longrightarrow (1)$: Let $\NN$ be a humble extension of $\MM$ such that $\Lt(\NN / \MM) \cong {\mathbf 2} \oplus L$. Let $b$ generate $\NN$ over $\MM$ such that $b$ fills the standard cut.
 
For every nonstandard $c \in M$, there is a nonstandard $n \in M$ such that, in $\MM$, there is an $n$-CPP representation $\alpha : L \into \eq(A)$, where $|A| \leq c$. Thus, for every $n < \omega$ and every nonstandard $c \in M$,
there is an $n$-CPP representation $\alpha : L \into \eq(A)$, where $|A| \leq c$. By underspill, for every $n < \omega$, there is, in $\MM$, $n$-CPP representation $\alpha : L \into \eq(A)$, where $|A|$ is standard. But such an $\alpha$ actually is an $n$-CPP representation of $L$. Therefore, $(1)$ holds.

\smallskip

 $(1) \Longrightarrow (2)$:  By Lemma~2.8, let $\NN_0 \succ \MM$ be a minimal, humble extension of $\MM$ that is generated over $\MM$ by $a$ filling  the standard cut of $\MM$. Since $L$ has $n$-CPP representations for each $n < \omega$, then $L$ has finite $n$-CPP representations for each $n < \omega$. By overspill, there is a nonstandard $k \in N_0$ and, in the sense of $\NN_0$, an $\NN_0$-finite $k$-CPP representation $\alpha : L \into \eq(A)$. Moreover, we can arrange $k$ and $\max(A)$ fill the standard cut of $\MM$. Now modify 
 $\alpha$ to get $\beta : L \into \eq(B)$, where $B = \{a\} \times A$ and $\alpha \cong \beta$ by the bijection that maps $x$ to $\langle a,x \rangle$. Thus, $\beta$ is $k$-CPP. Also, there  function $t : B \into N_0$ that is constantly $a$ is $\MM$-definable. We can now get $\NN \succ \NN_0$ such that $\Lt(\NN / \NN_0) \cong L$ and there is $b$ that generates $\NN$ over $\NN_0$ such that $b \in B^\NN$.  Then $\NN \models t(b) = a$, so $\Lt(\NN / \MM) \cong {\mathbf 2} \oplus L$. \qed

 \bigskip

 The previous lemma suggests the question: If $\NN$ is a humble extension of $\MM$ and $\Lt(\NN / \MM)$ is finite, is there a lattice $L$ such that $\Lt(\NN / \MM) \cong {\mathbf 2} \oplus L$? The answer is {\it not always}; in fact, if $\MM$ is nonstandard and countable and $({\mathbb N}, \ssy(\MM) )\models {\sf ACA_0}$, then there is a humble extension $\NN \succ \MM$ such that $\Lt(\NN / \MM) \cong {\mathbf B}_2$, where  ${\mathbf B}_2$ is the Boolean lattice having exactly 2 atoms. We sketch a proof that uses   the Canonical Ramsey Theorem for pairs (${\sf CRT}^2$), which is  a consequence of Ramsey's Theorem for 4-tuples  \cite{er}, which, in turn, is  well known  to be a consequence of ${\sf ACA}_0$. (Mileti \cite{m} has proved the converse: ${\sf RCA}_0$ + ${\sf CRT}^2$ implies ${\sf ACA}_0$.)
 
 \begin{quote} 
 ${\sf CRT}^2$: If  $f : \omega^2 \into \omega$, then there is an unbounded $X \subseteq \omega$ such that one of the following:
 
 \begin{itemize}
 
 \item $f$ is either constant  or one-to-one on $X^2$;
 
 \item there is a one-to-one function $g : \omega \into \omega$ such that either
 
 \begin{itemize} 
 
 \item whenever $x,y \in X$ and $x < y$, then $f(x,y) = g(x)$, or
 
  \item whenever $x,y \in X$ and $x < y$, then $f(x,y) = g(y)$.
  
  \end{itemize}
  \end{itemize}
  \end{quote}
  Let $X_0 \supseteq  X_1 \supseteq X_2 \supseteq \cdots$ be a sequence of sets in $\Def(\MM)$ such that whenever $f : M^2 \into M$ is $\MM$-definable such that $f(n) < \omega$ whenever $n < \omega$, then there is $i < \omega$ such that if $X = X_i \cap \omega$, then $X$ is as in  the conclusion of  ${\sf CRT}^2$. This sequence determines a humble extension $\NN \succ \MM$ for which $\Lt(\NN / \MM) \cong {\mathbf B}_2$. 
  
  \bigskip
  
  {\sc Question 2.10}: What conditions on $\MM$ (if there are any) guarantee that whenever $\NN$ is a humble extension of $\MM$ and $\Lt(\NN / \MM)$ is finite, then 
  $\Lt(\NN / \MM) \cong {\mathbf 2} \oplus L$ for some $L$?

 \bigskip

 Returning to the main thread, we next consider ordertypes. Two linearly ordered sets $(A,<)$ and $(B,<)$ are isomorphic iff they have the same ordertypes. A finite linearly ordered set $(A,<)$ with $|A| = n$ has ordertype $n$. We let ${\mathbb Z}$ be the ordertype of the (positive and negative) integers and $\omega$ be the ordertype of the nonnegative integers. 
 If $\alpha$ and $\beta$ are the ordertypes of $(A,<)$ and $(B,<)$, respectively, then $\alpha + \beta$ is the ordertype of the disjoint union of $A$ and $B$, where $a < b$ whenever $a \in A$ and $b \in B$. If $\lambda_0,\lambda_1, \lambda_2, \ldots$ are order types, then we define the ordertype 
 $\lambda_0 + \lambda_1 + \lambda_2 + \cdots$ in the obvious way. Every linearly ordered set can be thought of as a lattice. 
 
 Suppose that $L$ is a finite lattice. We make use of $n$-CPP representations $\alpha : L \into \eq(A)$.
The ${\mathcal L}$-formula $cpp_L(x)$ is such that if $a \in M$, then $\MM \models cpp_L(a)$ iff $\MM$ thinks that $L$ has  an $a$-CPP representation. 
 
  \bigskip

{\sc Theorem 2.11}: {\em Let $\MM$ be countable and nonstandard, and let $\lambda$ be a countable ordertype $($allowing the possibility that $\lambda = 0)$. Then $\MM$ has a humble extension $\NN$ such that 
$$ 
\Lt_0(\NN / \MM) \cong 2+ \lambda + 1.
$$} 

\bigskip

{\it Proof}. Let $L$ be a linearly ordered set having ordertype $2+\lambda + 1$. Let $L = \{a_i : i < \omega\}$, where $a_0$ is the first element, $a_1$ is the second, and $a_2$  the last element. The  $a_i$'s need not be distinct, which will definitely be the case if $L$ is finite; 
it is even possible that   $L =  \{a_0,a_1,a_2 \}$. 

Just as in the proof of $(1) \Longrightarrow (2)$ of Lemma~2.9, we let $\NN_0 \succ \MM$ be a minimal, humble extension that is generated  over $\MM$ by $a$ that fills the standard cut of $\MM$. Let ${\mathbf L}$ be the class of all finite linearly ordered sets, thought of as lattices. If $L \in {\mathbf L}$, then an $\NN_0$-definable representation $\alpha : L \into \eq(A)$  is {\it ample} if whenever $r < s \in L$, then every $\alpha(r)$-class includes a nonstandard number of $\alpha(s)$-classes. Recall that each lattice in ${\mathbf L}$ has finite $n$-CPP representations for each $n < \omega$. Thus, if $L \in {\mathbf L}$ and $|L| \geq 2$, then $\NN_0$ thinks that each ample representation is $n$-CPP for each $n < \omega$.

For $i < \omega$, let $L_i = \{a_1, \ldots, a_{i+2}\} \subseteq L$.  Thus, each $L_i \in {\mathbf L}$ and $L_0$ has ordertype  $2$. Let $\alpha_0 : L_0 \into \eq(A_0)$ be an ample representation of $L_0$,
where $A_0 = \{a\} \times A$ for some $\NN_0$-bounded $A \in \Def(\NN_0)$.  Let 
$N_0 \backslash \omega = \{d_0,d_1, \ldots \}$, and let $\Theta_0,\Theta_1,\Theta_2, \ldots$ enumerate all the $\NN_0$-definable $\Theta \in \eq(M)$. 
For each $i < \omega$, let $\alpha_{i+1} : L_i \into \eq(A_{i+1})$ be an ample representation  of $L_{i+1}$ such that:

\begin{itemize}

\item $A_{i+1} \subseteq A_i$ and $\alpha_{i+1} \harp L_i = \alpha_i | A_i$,

\item $|A_{i+1}| \leq d_i$,

\item $\Theta_i$ is canonical for $\alpha_{i+1}$.

\end{itemize} 
The sequence $A_0 \supseteq A_1 \supseteq A_2 \supseteq \cdots$ determines a complete 1-type $p(x)$ over $\NN_0$. Let $\NN$ be generated over $\NN_0$ by $c$ that realizes $p(x)$. Then $c$ fills the standard cut of $\MM$ and one easily verifies that $\Lt_0(\NN / \NN_0) \cong 1 + \lambda + 1$ and that $\Lt_0(\NN / \MM) \cong L$. \qed 

\bigskip

The following corollary to Theorem~2.11 will be used   in the proof of Theorem~B in \S3.

\bigskip

{\sc Corollary 2.12}: {\em  There is an ${\mathcal L}^*$-formula $\theta(x)$ such that whenever  $I \subseteq \omega$ and $\MM$ is countable and nonstandard,  then $\MM$ has a humble  extension $\NN \succ \MM$  such that   $\theta(x)$ defines $I$ in $(\NN,M)$.}
  
  \bigskip

   {\it Proof}. If $\NN$ is a humble extension of $\MM$, then $\omega$ is definable in $(\NN,M)$ since $\omega = \gcis(\NN,\MM)$. In fact, there is a single definition that works whenever $\NN$ is a humble extension of $\MM$. 
   
   Also, if $\NN$ is a humble extension of $\MM$, and $b \in N \backslash M$ generates $\MM(b)$ over $\MM$, then $M(b)$ is definable in $(\NN,M)$. Furthermore, the function $b \mapsto M(b)$ is definable in $(\NN,M)$. In fact, there is a single definition that works whenever $\NN$ is a humble extension of $\MM$.

   For  $I \subseteq \omega$, let 
$$
\lambda_I =  \omega + \big((n_0 + {\mathbb Z}) + (n_1 + {\mathbb Z}) + (n_2 + {\mathbb Z})  + \cdots \big)+ 1,
$$ 
where $n_i = 1$ if $i \in I$ and $i = 2$ if $i \not\in I$. Thus, $\lambda_I$ is a countable ordertype with  first, second and last elements, so that by Theorem~2.11, there is a humble extension $\NN \succ \MM$ such that 
$$
\Lt_0(\NN / \MM) \cong \lambda_I.
$$

Now, it is easily seen that there is a formula $\theta(x)$ that defines $I$ in $(\NN,M)$ whenever $\NN$ is a humble extension of $\MM$ and  $\Lt_0(\NN / \MM) \cong \lambda_I$. This $\theta(x)$ is as desired. \qed

\bigskip

An example of an ordertype not covered by Theorem~2.11 is $1 + {\mathbb Z} + 1$. Nevertheless, it is true that every countable nonstandard $\MM$ has a humble extension $\NN$ such that $\Lt_0(\NN / \MM) \cong 1 + {\mathbb Z} + 1$. The following improvement of Theorem~2.11 has a proof that is reminiscent of the proof of \cite[Th.~3.2/Coro.~3.3]{almost}. Since we can get away without this improvement, we present only an abbreviated proof.

\bigskip

{\sc Theorem 2.13}: {\em Let $\MM$ be countable and nonstandard, and let $\lambda$ be a countable ordertype $($allowing the possibility that $\lambda = 0)$. Then $\MM$ has a humble extension $\NN$ such that 
$$ 
\Lt_0(\NN / \MM) \cong 1+ \lambda + 1.
$$} 

\bigskip

{\it Proof}. Let $L$ be a linearly ordered set having ordertype $1 + \lambda +1$. Let $L = \{a_i : i < \omega\}$, where $a_0$ is the first element and $a_1$ is the last. For $i < \omega$, let 
$L_i = \{a_0,a_1, \ldots, a_{i+1}\}$. Let's say,  for $i < \omega$, that a representation $\alpha : L_i \into \eq(A)$ is {\it ample} if the following:

\begin{itemize}

\item $\alpha$ is $\MM$-definable (so $A \in \Def(\MM)$).

\item $A \cap \omega$ is infinite.

\item  Whenever $b \in A \cap \omega$ and $1 \leq j \leq i$ there are $a,c \in A \cap \omega$ such that 
$a \leq b < c$ and $[a,c) \cap A$ is an $\alpha(a_j)$-class.

\item For each $n < \omega$, there is an $\alpha(a_1)$-class $A \subseteq \omega$ such that whenever $1 \leq j \leq i$ and $B \subseteq A$ is an $\alpha(a_j)$-class, then $B$ includes at least $n$ 
$\alpha(a_{j+1})$-classes.

\end{itemize} 
There are three relevant facts about ample representations that will be needed.

\smallskip

{\sf Fact 0}: There is an ample representation of $L_0$.

\smallskip

{\sf Fact 1}: Whenever $\alpha : L_i \into \eq(A)$ is ample and $\Theta \in \eq(A)$ is $\MM$-definable,
then there is $B \subseteq A$ such that $\alpha|B$ is ample and $\Theta$ is canonical for $\alpha|B$.

\smallskip

{\sf Fact 2}: Whenever $\alpha : L_i \into \eq(A)$ is ample then there is an ample $\beta : L_{i+1} \into \eq(A)$ such that  $\alpha = \beta \harp L_i$.

\smallskip
Making use of the countability of $\MM$, we do the usual construction of  a decreasing sequence $X_0 \supseteq X_1 \supseteq X_2 \supseteq \cdots$ of sets in $\Def(\MM)$ such that for each $i < \omega$, there is an ample representation $\alpha_i : L_i \into \eq(X_i)$ such that: 
\begin{itemize}

\item For each $i < \omega$, 
$ \alpha_i | X_{i+1} = \alpha_{i+1} \harp L_i$.

\item For every $\MM$-definable $\Theta \in \eq(M)$, there is $i < \omega$ such that $\Theta$ is 
canonical for $\alpha_i$. 

\end{itemize}
There is a unique ultrafilter ${\mathcal U}$ of $\ssy(\MM)$ such that $X_i \in {\mathcal U}$ for each $i < \omega$. Then $\NN = \MM^\omega/{\mathcal U}$ is as required. \qed

%
%

{\section{Some Combinatorics}} Only  finite combinatorics is discussed in this section; other than in the first and  last paragraphs, there is no mention of \pa\ at all. 

We begin by confirming some standard notation. Each nonnegative integer $n$ is such that  $n = \{0,1,2, \ldots, n-1\}$. As usual, $\omega$ is the set of nonnegative integers, and we  usually write $n < \omega$ instead of $n \in \omega$. For any linearly ordered set $A$,  we will use  the interval notation for  $k \leq \ell$: $[k,\ell] = \{m \in A : k \leq m \leq \ell\}$ and $[k,\ell) = \{ m \in A : k \leq  m< \ell\}$. In particular, $n = [0,n)$ for  $n < \omega$, although for such an $n$, we prefer $n$ to $[0,n)$. If $n < \omega$ and $A$ is a set, then $A^n$ is the set of $n$-tuples of elements of $A$ (except that if $k < \omega$, then $k^n$ is given its usual arithmetical meaning).  Equivalently, $A^n$ is the set of functions $f : n \into A$ (so $A^0 = \{\varnothing\}$ if $A \neq \varnothing$).   
Note especially that  $\omega^n$ is a set of $n$-tuples and not an ordinal.

We now give the definitions of  those objects that will occupy our attention for the remainder of this section.

\bigskip

{\sc Definition} 3.1a:  For $1 \leq h < \omega$, we  say that $(X;  E_0,E_1, \ldots,E_{h-1})$ is an $h$-{\bf structure} if $X$ is a nonempty finite set,  $E_i \in \eq(X)$ for each $i < h$, $E_0$ is trivial, and whenever $i < h-1$, then 
$E_i \supseteq E_{i+1}$. If $(X;E_0,E_1, \ldots,E_{h-1})$ is an $h$-structure and $\varnothing \neq Y \subseteq X$, then $Y$ \mbox{{\bf induces}} the 
$h$-substructure $(Y; E_0 \cap Y^2,E_1 \cap Y^2, \ldots,E_{h-1} \cap Y^2)$. 

\bigskip

We often use $(X; \overline E)$ as shorthand for the $h$-structure $(X;  E_0,E_1, \ldots,$ $E_{h-1})$.
If $(X;\overline E)$ is an $h$-structure, then $E_i$ is defined only for $i < h$; nevertheless, we will sometimes refer to $E_h$, which should be  understood to be the discrete equivalence relation on $X$.

\bigskip

{\sc Definition 3.1}b: We will say that the $h$-structure  $(X; \overline E)$ is $(I,w)$-{\bf special} if $w \geq 2$ and the following:

\begin{itemize}

\item [(1)]  $I = \big\{ i < h : E_i$ is not discrete and $i = \max(\{j < h : E_i = E_j \})\big\}$;

\item[(2)] if $i \in I$ and  $A$ is an $E_i$-class, then 
\begin{equation*}
|\{B \subseteq A : B {\mbox{ is an }} E_{i+1}{\mbox{-class}}\}| \geq w. 
\end{equation*}

\end{itemize}

\bigskip

There is a close connection between $h$-structures and representations of the lattice ${\mathbf n}$, where $n = h+1$. (The lattice ${\mathbf n}$ is defined in the introduction. The elements of ${\mathbf n}$, in order, are $0_{\mathbf n} = 0 < 1 < \cdots < n-1 = 1_{\mathbf n}$.) Given a finite representation $\alpha : {\mathbf n} \into \eq(A)$, then $(A;{\overline E})$ is an $h$-structure, where $E_i = \alpha(i)$ for $i < k$. Conversely, if $(A;{\overline E})$ is an $(h,2)$-special $h$-structure, then $\alpha : {\mathbf n} \into \eq(A)$ is a representation of ${\mathbf n}$, where $n = h+1$, $\alpha(i) = E_i$ for $i < h$ and $\alpha(h)$ is discrete. Since we know that ${\mathbf n}$ has $k$-CPP representations for each $k < \omega$, we can conclude that for each $k < \omega$, there is $w < \omega$ such that if $(A;\overline E)$ is an $(h,w)$-special $h$-structure, then the corresponding 
$\alpha : {\mathbf n} \into \eq(A)$ is a $k$-CPP representation of ${\mathbf n}$. One goal of this section is get a useful upper bound on $w$ in terms of $h$ and $k$.

\bigskip

If $(X;\overline E)$ is an $(I,w)$-special $h$-structure, then $|X| \geq w^{|I|}$. For certain purposes, an  $(I,w)$-special $h$-structure is sufficiently similar to a $(k,w)$-special $k$-structure, where $k = |I|$.
We will often take advantage of this similarity without explicitly mentioning it. 

\bigskip

Obviously, there are plenty of $h$-structures. If $2 \leq h < \omega$, then let $X$ be the set of all functions $f : h \into h$. For each $i < h$, let $E_i \in \eq(X)$ be such that if $f,g \in X$, then 
$\langle f,g \rangle \in E_i$ iff $f \harp i = g \harp i$. Clearly, $(X; \overline E)$ is an $h$-structure; it is  even an $(h,h)$-special $h$-structure. We will call it the {\bf basic} $h$-structure.

Each part of the next proposition has a very routine proof, which is omitted.  This proposition is included here for easy future reference.

\bigskip

{\sc Proposition} 3.2:  {\em Let $(X; \overline E)$ be an $(I,w)$-special $h$-structure. 

$(1)$ If $\varnothing \neq I'\subseteq I$ and $2 \leq v \leq w$, then there is $Y \subseteq X$ that induces an $(I',v)$-special $h$-substructure such that $|Y| = v^{|I'|}$.

$(2)$   If $h \leq h' < \omega$, $I \subseteq I' \subseteq h'$ and $w \leq w' < \omega$, then there is an  $(I',w')$-special $h'$-structure 
$(X'; \overline{ E'})$ such that $X' \supseteq X$  and $E_i = E'_i \cap X^2$ for $i \leq h$.} \qed

\bigskip

{\sc Lemma} 3.3: {\em Suppose that $2 \leq w < \omega$ and $1 \leq c < \omega$. If $(X;{\overline E})$ is an $(I,cw)$-special $h$-structure and $f : X \into c$, then 
 there is $Y \subseteq X$ such that  $f$ is constant on $Y$ and $Y$ induces an $(I,w)$-special $h$-substructure of $(X; {\overline E})$.}

\bigskip

{\it Proof}. It suffices to assume that $I = h$ so that $(X;{\overline E})$ is an $(h,cw)$-special $h$-structure. We give a proof by induction for $h \geq 1$.

\smallskip

Basis step: $h = 1$. Let $(X;E_0)$ be a $(1,cw)$-structure. Since $E_0$ is trivial, that means that $|X| \geq cw$.    Let $f : X \into c$.  Since $|X| \geq cw$,  {\sf PHP1} implies there  is $Y \subseteq X$ such 
that $|Y| = w$ and  $f$ is constant  on $Y$. Then, $Y$ induces a $(1,w)$-substructure of  $(X;E_0)$.

\smallskip

Inductive step: Let $1 \leq h < \omega$, $2 \leq w < \omega$, $(X; \overline E)$ be an \mbox{$(h+1,cw)$}-special $(h+1)$-structure and  $f : X \into c$. For each $E_1$-class $B$, $B$ induces an $((h+1) \backslash \{0\}, cw)$-special $(h+1)$-substructure. 
 By the inductive hypothesis, let $Y_B \subseteq B$ and $d_B < c$ be such that $f$ is constantly $d_B$ on $Y_B$ and $Y_B$ induces an $((h+1) \backslash \{0\}, w)$-special $(h+1)$-substructure. Since there are at least $wc$ distinct $E_1$-classes,  {\sf PHP}1 implies that there are $d < c$ and a set ${\mathcal B}$ of $E_1$-classes such that  $|{\mathcal B}| = w$ and $d_B = d$ for each $B \in {\mathcal B}$.  Let 
 $Y = \bigcup\{Y_B : B \in {\mathcal B}\}$. Then, $f$ is constantly $d$ on $Y$, and $Y$ induces a $(k+1,w)$-special $(h+1)$-substructure of  $(X;{\overline E})$. \qed
 
 \bigskip

 {\sc Lemma 3.4}: {\em   Suppose that $(X; \overline E)$ is an $(I,w^{k+1})$-special $h$-structure, $k = \max(I)$ and  $f : X \into \omega$ is such that $f$ is one-to-one on each $E_k$-class.  Then there is $Y \subseteq X$ that induces an $(I, w)$-special $h$-substructure such that $f$ is one-to-one on $Y$.}
 
 \bigskip
 
 {\it Proof}. It suffices to assume that $I = h$, so that $k = h-1$.  Thus, we have that $(X;\overline E)$ is an  $(h,w^h)$-special $h$-structure. If $h = 1$, then $f $ is one-to-one on $X$, so let $Y = X$. Thus, we  assume that $h \geq 2$. 
 
We easily get $Z \subseteq X$ that induces an $h$-substructure so that: 
 \begin{itemize}
 
 \item [(1)] Each $(E_{h-1} \cap Z^2)$-class has exactly $w^h$ elements.
 
 \item[(2)] For each $i < h-2$,  each $(E_i \cap Z^2)$-class includes exactly $w$ $(E_{i+1} \cap Z^2)$-classes.
 
 \end{itemize} 
 
 From (2), there are exactly $w^{h-1}$ $(E_{h-1} \cap Z^2)$-classes. 
   Let $\{A_i : i < w^{h-1}\}$ be the set of $(E_{h-1} \cap Z^2)$-classes. Since $f$ is one-to-one on each $A_i$ and each $A_i$ has $w^h$ members, we  can obtain, by recursion on $i$, subsets $B_i \subseteq A_i$ such that 
 $|B_i| = w$ and $\{f(x) : x \in B_i\} \cap \{f(x) : x \in B_j,j<i\} = \varnothing$. Let $Y = \bigcup\{B_i : i < w^{h-1}\}$. Then $Y$ induces a $(h,w)$-special $h$-substructure of $(X;\overline E)$ such that $f$ is one-to-one on~$Y$. \qed
 
 \bigskip

Let $(X; \overline E)$ be an $h$-structure.  A function $f$ whose  domain includes $X$ is  $f$  {\bf canonical} on $(X; \overline E)$ if   there is $i \leq h$  such that whenever $x,y \in X$, then $f(x) = f(y)$ iff $\langle x,y \rangle \in E_i$. If we wish to point out the $i$ in this definition, then we say that $f$ is $i$-{\bf canonical}.   Notice that if $(X;\overline E)$ is an $(I,w)$-special $h$-structure and $f$ is canonical on $(X;\overline E)$ then $f$ is $i$-canonical for some $i \in I \cup \{h\}$. 
Furthermore, $f$ is $h$-canonical iff it is one-to-one on $X$.

 \bigskip
 
 {\sc Lemma 3.5}: {\em  Suppose that $(X; \overline E)$ is an $(I,w^{k+1})$-special $h$-structure and $k \in I$. Let $f : X \into \omega$ be such that $f$ is constant on each $E_{k+1}$-class and whenever $C$ is an $E_k$-class,  $A,B \subseteq C$ are $E_{k+1}$-classes, $x \in A$ and $y \in B$, then $f(x) = f(y)$ iff $A = B$.  
  Then there is $Y \subseteq X$ that induces an $(I, w)$-special $h$-substructure on which $f$ is $k$-canonical.}
  
   \bigskip
 
 {\it Proof}. If $k = \max(I)$, then this lemma reduces to Lemma~3.4. So, assume that $k < \max(I)$. It suffices to assume that $I = h$, so that $k \leq h-1$. Thus, $(X;\overline E)$ is an $(h,w^{k+1})$-special $h$-structure. 
 
 Let $Z$ be the set of $E_k$-classes, and then let $(Z; \overline F)$ be the $k$-structure such that if $i < k$ and $A,B \in Z$, then $\langle A,B \rangle \in F_i$ iff $\langle x,y \rangle \in E_i$ for all (or, equivalently,  some) $x \in A$ and $y \in B$. We easily see that $(Z; \overline F)$ is a $(k,w^{k+1})$-special $k$-structure.
 Let $g : Z \into \omega$ be such that if $A \in Z$, then $g(A) = f(x)$, where $x \in A$.
 Since $f$ is constant on each $A \in Z$, $g$ is well defined. Moreover,  $g$ is one-to-one on each 
 $F_k$-class. By Lemma~3.4, there is $Z_0 \subseteq Z$ that induces an $(k,w)$-special 
 $k$-substructure of $(Z; \overline F)$ such that $g$ is one-to-one on $Z_0$. It is easily verified that $Y = \bigcup Z_0$ induces an $(h,w)$-special $h$-substructure of $(X; \overline E)$ on which $f$ is $k$-canonical. \qed

 \bigskip
 
 If $ 1 \leq h < \omega$, then the function $C_h : \omega \into \omega$ is defined so that \big($C_h(0) = 0$ and $C_h(1) = 1$, to make it total,  and\big) whenever $2 \leq w < \omega$, then $C_h(w)$ is the least (if it exists) such that $2 \leq C_h(w) < \omega$ and whenever $(X;{\overline E})$ is an $(h,C_h(w))$-special $h$-structure and $f : X \into \omega$, then there is $Y \subseteq X$ that   induces an $(h,w)$-special $h$-substructure on which $f$ is canonical. If $w < v < \omega$, then every $(h,v)$-special $h$-structure is $(h,w)$-special; therefore, if $w < v < \omega$ and $C_h(v)$ exists, then so does $C_h(w)$ and $C_h(w) \leq C_h(v)$. 
 
 The next lemma shows that $C_h(w)$ always exists and gives upper bounds.

\bigskip

{\sc Lemma} 3.6: {\em If $2 \leq w < \omega$, then$:$

$(a)$ $C_1(w) \leq w^2;$

$(b)$ if $1 \leq h < \omega$, then $C_{h+1}(w) \leq C_h((h+1)w^{h+1})$.}

\bigskip

{\it Proof}. (a): Let $(X; E_0)$ be a $(1,w^2)$-special $1$-structure. Since $E_0$ is trivial, then  $|X| \geq w^2$. Let $f : X \into \omega$. Then by {\sf PHP2}, there is $Y \subseteq X$ such that $|Y| = w$ (so that $Y$ induces a $(1,w)$-special $1$-substructure) and $f$ is either constant on $Y$ (i.e. $0$-canonical) or $f$ is one-to-one on $Y$ (i.e. $1$-canonical).

\smallskip

(b): Let $1 \leq h < \omega$. Let $(X; \overline E)$ be an $\big(h+1,C_h((h+1)w^{h+1})\big)$-special $(h+1)$-structure, and let $f : X \into \omega$. Each $E_1$-class $A$ induces a  $\big([1,h],C_h((h+1)w^{h+1})\big)$-special $(h+1)$-substructure of $(X;\overline E)$. For each such $A$, let $B_A \subseteq A$ induce a $([1,h], (h+1)w^{h+1})$-special $(h+1)$-substructure on which $f$ is $j_A$-canonical, where $1 \leq j_A \leq  h+1$.  By {\sf PHP1}, let $1 \leq j \leq h+1$ and ${\mathcal B}$ be a set of $E_1$-classes such that 
$|{\mathcal B}| = w^{h+1}$ and $j_A = j$ for each $A \in \mathcal B$. 
There are two cases, depending on whether $j = 1$ or $2 \leq j \leq h+1$.

\smallskip

{\sf Case 1}: $j=1$. Thus, for each $A \in {\mathcal B}$, $f$ is constant on $B_A$. Let $g : {\mathcal B} \into \omega$ be such that 
$g(A) = f(x)$ for some (or all) $x \in B_A$. By {\sf PHP2}, let ${\mathcal{B'}} \subseteq \mathcal B$ be such that $|{\mathcal{B'}}| = w$ and $g$ is either constant or one-to-one on ${\mathcal{B'}}$. Let $Y = \bigcup \{B_A : A \in \mathcal{B'}\}$. Then $Y$ induces   an $(h+1,w)$-special $(h+1)$-substructure.
Then $f$ is $0$-canonical on this $(h+1)$-substructure (if $g$ is constant on $\mathcal{B'}$) or $f$ is $1$-canonical on it 
(if $g$ is one-to-one on $\mathcal{B'}$).

\smallskip

{\sf Case 2}: $2 \leq j \leq h+1$. Let $Z = \bigcup {\mathcal B}$. Thus, $Z$ induces an 
$(h+1,w^{h+1})$-special $(h+1)$-substructure $(Z; \overline F)$ of $(X;\overline E)$. Whenever $\langle x,y \rangle \in  F_{j-1}$, then $f(x) = f(y)$ iff $\langle x,y \rangle \in F_j$. Lemma~3.5 implies that that there is $Y \subseteq Z$ that induces  
an $(h+1,w)$-special $(h+1)$-substructure on which $f$ is $j$-canonical. 
 \qed
 
 \bigskip 
 
 The following corollary gives an upper bound for $C_h(w)$ which is good enough for our purposes.
 
 \bigskip

{\sc Corollary} 3.7: {\em If $1 \leq h <  \omega$ and $2 \leq w < \omega$, then 
$$
C_h(w) \leq (h!w)^{2h!}
$$.}

\smallskip

{\it Proof}. This is a routine consequence of Lemma~3.6. \qed

\bigskip

{\sc Corollary} 3.8: {\em Suppose that $1 \leq k < \omega$, $2 \leq w < \omega$ and  $(X; \overline E)$ is an $(I,(k!w)^{2k!})$-special $h$-structure, where $k = |I|$ and $f : X \into \omega$. Then there is  $j \in I$ such that$:$

\begin{itemize}

\item If $j \neq \max(I)$, then there is $Y \subseteq X$ on which $f$ is constant and $Y$ induces an $(I \cap [j,h),w)$-special $h$-substructure. 

\item If $j \neq \min(I)$, then there is $Z \subseteq X$ on which $f$ is one-to-one and $Z$ induces an $(I \cap [0,j),w)$-special $h$-substructure.  
\end{itemize}}

\bigskip

{\it Proof}. Suppose that $(X; \overline E)$ is an $(I,(k!w)^{2k!}$-special $h$-structure, where $k = |I|$.  Let  $f : X \into \omega$. Notice that the corollary is vacuously true if $k = 1$. So suppose that $k \geq 2$. By Corollary~3.7, $(k!w)^{2k!} \geq C_k(w)$, so we can let $X_0 \subseteq X$  induce an $(I,w)$-substructure of $(X; \overline E)$ on which $f$ is canonical. Thus, there is  $j' \in I \cup \{k\}$ such that  $f$ is $j'$-canonical. 

If $j' = k$, then $f$ is one-to-one on $X_0$, so we can we can let $j = \max(I)$ and $Z = Z_0$.

If $j' \neq k$, then let $j = j'$. By Proposition~3.2(2), let $Y \subseteq Z_0$ induce an $(I \cap j,h),w)$-special $h$-substructure, and if $j \neq \min(I)$, let $Z \subseteq Z_0$ induce an $(I \cap [0,j))$-special $h$-substructure. These $Y$ and $Z$ clearly are as required.  \qed

\bigskip

Each definition and result in  this section can be formalized in \pa. Furthermore,  the proof of each of the results is straightforwardly transferred to a \pa-proof  of its formalization. We will be making use of the formalizations of some of these results in the next sections. As a typical example of how they  will be referred to, the formalized version of    Lemma~3.5 will be  Lemma~3.5(\pa).

\bigskip
%
%

\section{Proving Theorem~B}  Theorem~B will be proved in this section.  We begin with some general remarks about   minimal extensions as in Theorem~B while, for the time being, ignoring the ${\mathcal L}^*$-formula $\varphi(x)$. 

In the proofs of $(1) \Longrightarrow (3)$ of Theorems~1.1 and~1.2, we saw how to get an extension 
$\NN \succ \MM$ such that (M1) -- (M4). In our proof of Theorem~B, we will do a similar construction but being more careful about the  $X_i$'s. 

Suppose that $J$ is an exponentially closed cut of $\MM$. We work in $\MM$ when it is appropriate. 

For each $2 \leq h \in M$, we let $P_h \in \Def(\MM)$ be  (the code of) the set of $\MM$-definable functions $f : h \into h$. Obviously, $P_k \cap P_h = \varnothing$ for distinct $k,h$. Recall that for each such $P_h$, there is a basic $h$-structure $(P_h;\overline E)$. In this section, whenever we are considering some nonempty $Y \subseteq P_h$, we will often be identifying it with the $h$-substructure it induces. We will say that  $Y$ is {\bf viable} if $Y \in \Def(\MM)$, $Y \subseteq P_h$ and $Y$ induces an $(I,w)$-special $h$-substructure
  such that 
$J < w \leq h$ and $I \cap \omega$ is infinite. 
 For such a viable $Y$, we sometimes say that $Y$ is $I$-{\bf viable} or $(I,w)$-{\bf viable} if we wish to indicate $I$ and/or $w$. 

It is obvious that $P_h$ is viable for $h \in M \backslash J$; in fact, it is    $([0,h),h)$-viable.    Every viable $Y$ is large as defined in the proof of Theorem~1.1.

\bigskip

{\sc Lemma 4.1}: {\em Suppose that $J$ is an exponentially closed cut of $\MM$ and that $Y$ is viable.

\begin{itemize} 

\item[(V1)] If  $f : Y \into M$ is $\MM$-definable, then there is a viable $Z \subseteq Y$ such that $f$ is either constant or one-to-one on $Z$.

\item[(V2)] If $d \in M \backslash J$, there is a viable $Z \subseteq Y$ such that $\MM \models |Z| \leq d$.

\item[(V3)] If  $f : Y \into M$ is  one-to-one and $\MM$-definable, there are a viable $Z \subseteq Y$ and $d \in M \backslash J$ such that $\MM \models \forall x \in Z [f(x) \geq d]$.

\item[(V4)] For every $\MM$-definable function $f : Y \into \Def(\MM)$, there are viable $Z \subseteq Y$ and $d \in M \backslash J$ such that for all $x,y \in Z$, $f(x) \cap [0,d) = f(y) \cap [0,d)$.

\end{itemize}
Furthermore, if $Y$ is $I$-viable, then in each of \em{(V1) -- (V4)}, {\em there is such a $Z$ that is $K$-viable, where    $K \cap (\omega \backslash m) = I \cap (\omega \backslash m)$  for some $m < \omega$.} 

\bigskip

{\it Proof}. Without loss of generality,  assume that $Y \subseteq P_h$ and that $Y$ is $([0,k),w)$-viable. Thus, 
$k,w \in M$ are such that $k$ is nonstandard and $w > J$.
\smallskip

(V1): Let $f: Y \into M$ be $\MM$-definable. Since $J$ is an exponentially closed cut and $w > J$, there is $v > J$ such that $v^v < w$. For each $\ell < \omega$, we have that 
\begin{equation}
v \leq C_\ell(v) \leq v^v < w. \tag{$*$}
\end{equation}
 The first inequality is trivial, and the second follows easily from Corollary~3.7(\pa). By overspill, there is a nonstandard $\ell > \omega$ for which $(*)$ holds. Without loss (by Proposition~3.2(\pa)), we can assume that $C_k(v) < w$. 

 Let $j < k$ be as in Corollary~3.8(\pa).  If $j$ is standard (so $j <k-1$), then let $Z \subseteq Y$ be such that $f$ is constant on $Z$ and $Z$ induces a $([j,k),w)$-special $h$-structure; if  $j$ is nonstandard (so $j \neq 0)$, then let $Z \subseteq Y$ be such that $f$ is one-to-one on $Z$ and $Z$ induces an $([0,j),w)$-special $h$-structure. Either way,    $Z$ is viable and as required.  
\smallskip

(V2): Let $d \in M \backslash J$. Let $v \in M$ be such that $v > J$ and $v^v <\min(d,w)$. Let nonstandard $\ell$ be such that 
$\ell < \min(v,k)$.  By Proposition~3.2(\pa), let $Z \subseteq Y$ induce a $([0,\ell),v)$-special 
$h$-substructure such that $|Z| = \ell^v$.  Thus, $Z$ is $([0,\ell),v)$-viable. Clearly, $\ell^v < v^v \leq d$.  Then $Z$ is viable and $\MM \models |Z| \leq d$.

 \smallskip

(V3): Let $f : Y \into M$ be  one-to-one  and $\MM$-definable.  Let $d \in M$ be such that $d > J$ and $(d+1)^2 < w$.  
Let $f_d : Y \into [0,d]$ be such that $f_d(x) = \min(f(x), d)$. Since $f_d$ is $\MM$-definable, there is, by Lemma~3.3(\pa), a subset $Z \subseteq Y$  such that $f_d$ is constant on $Z$ and $Z$ induces a $([0,k), d+1)$-special $h$-substructure. Thus, $Z$ is $([0,k), d+1)$-viable. 
Clearly, $f_d(x) = d$ for all $x \in Z$, so that $f(x) \geq d$ for all $x \in Z$.   Hence, $Z$ and $d$ are as required.

\smallskip

(V4): Let $f : Y \into \Def(\MM)$ be $\MM$-definable. Let $d \in M$ be such that $d > J$ and 
$4^d <w$.  Let $f_d$ be the function on $Y$ such that $f_d(x) = f(x) \cap [0,d)$. Thus, we can think of $f_d : M \into 2^d$. Since $f_d$ is $\MM$-definable, there is, by Lemma~3.3(\pa), a subset $Z \subseteq Y$  such that $f_d$ is constant on $Z$ and $Z$ induces a $([0,k), \lfloor w/2^d \rfloor)$-special $h$-substructure. Clearly, $\lfloor w/2^d \rfloor \geq 2^d-1 > J$, so that $Z$ is viable.

\smallskip

In the proofs of each of (V1) -- (V4), it is clear that the $Z$ satisfies the ``furthermore'' part of the lemma. \qed

\bigskip

Each of (V1) -- (V4) of the previous lemma corresponds, respectively, to (S1) -- (S4) in the proofs of Theorems~1.1 and~1.2. Therefore, if $\MM$ is countable and $J$ is an exponentially closed cut, then  we can get a sequence $X_0 \supseteq X_1 \supseteq X_2 \supseteq \cdots$ of viable subsets of $M$ satisfying (S1) -- (S4), resulting in $\NN \succ \MM$ satisfying (M1) -- (M4). 
We say that an extension $\NN$ of $\MM$ obtained in this manner as a {\bf V-extension} of $\MM$. To be more precise, $\NN$ is a V-extension of $\MM$ iff there is a sequence $X_0 \supseteq X_1 \supseteq X_2 \supseteq \cdots$ of viable sets in $\Def(\MM)$ satisfying (S1) -- (S4) such that $\NN$ is generated over $\MM$ by an element in $\bigcap_{i<\omega}X_i^\NN$. We say that $a \in N$ is a {\bf V-point} if it is such a generator of $\NN$ over $\MM$.  
   \bigskip
   
    {\sc Lemma} 4.2: {\em There is an ${\mathcal L}^*$-formula $\Omega(x)$ such that whenever $\NN$ is a V-extension of $\MM$, then $\Omega(x)$ defines $\omega$ in $(\NN,M)$.}

\bigskip

{\it Proof}. Let $X_0 \supseteq X_1 \supseteq X_2 \supseteq \cdots$ be a sequence of  viable sets satisfying (S1) -- (S4). Let $c \in N$ generate $\NN$ over $\MM$, where $c \in X_i^\NN$ for each $i < \omega$. Thus, $c$ is a V-point. Following Proposition~2.2(\pa), there are, for each nonstandard $a \in M$ a nonstandard $b< a$ and  $Y \in \Def(\MM)$ such that $Y$ is $[0,b)$-viable and $c \in Y^\NN$. On the other hand, for each standard $b \in M$, there is no $Y \in \Def(\MM)$ such that, for some $w \in M$,  $Y$ is 
a $([0,b),w)$-special $h$-substructure and $c \in Y^\NN$.  Therefore, $\omega$ is definable from $c$ and $h$ in $(\NN,M)$.  As is well known, this implies that $\omega$ is definable in $(\NN,M)$ without parameters. 

 It should be clear that the defining formula we have just obtained can be chosen independently of $\MM$, $\NN$,$J$ and $h$. That is, there is a single ${\mathcal L}^*$-formula $\Omega(x)$ that works  for any nonstandard $\MM$, any exponentially closed cut $J$ of $\MM$ and any V-extension $\NN$. Thus, for all such $\NN \succ \MM$,
$$
(\NN,M) \models \Omega(a) {\mbox{ iff }} a \in \omega
$$
for all $a \in N$. 

It easily follows that the set of V-points in N is definable in $(\NN,M)$, and that such a definition is independent of $\MM$, $J$, $\NN$ and $h$. \qed

\bigskip

Next, we have an  easy lemma, adding a couple of  conditions to those in Lemma~4.1.
   
   \bigskip
   
   {\sc Lemma 4.3}: {\em Suppose that $J$ is an exponentially closed cut of $\MM$ and that $Y$ is  $I$-viable.

   \begin{itemize}
   
     \item[(1)]     If $K \in \Def(\MM)$ and $I \subseteq K \subseteq [0,h)$, then there is a $K$-viable $Z \supseteq Y$.  
   
   \item[(2)] If $K \in \Def(\MM)$, $K \subseteq I$ and $K \cap \omega$ is infinite, then there is a $K$-viable $Z \subseteq Y$.

   \end{itemize}}
   
   \bigskip
   
   {\it Proof}. This follows easily from Proposition~3.2(\pa). \qed
   
   \bigskip

   This lemma has  interesting consequences for which another definition is needed.    If $\NN$ is a V-extension of $\MM$ and $a \in N$ is a V-point,   define ${\mathcal F}_a$ to be the set of all $B \subseteq \omega$ for which there is an $I$-viable $Y$ such that $a \in Y^\NN$ and $B = I \cap \omega$. 
   
   \bigskip
   
   {\sc Corollary 4.4}: {\em Suppose that $J$ is an exponentially closed cut of~$\MM$. 
 \begin{itemize}
   \item[$(1)$] If $\NN$ is a V-extension of $\MM$ and $a \in N \backslash M$ is
     a V-point, then 
   ${\mathcal F}_a$ is a free filter of $\ssy(\MM)$.
   
  \item[$(2)$] If ${\mathcal F}$ is a free filter of $\ssy(\MM)$, then there is a V-extension $\NN$ of $\MM$ that is generated over $\MM$ by $a \in N$ such that ${\mathcal F}_a = {\mathcal F}$.\qed
  \end{itemize}}   
  
  \bigskip
  
  Part $(1)$ of this corollary shows that ${\mathcal F}_a$ is always a free filter (given the appropriate hypotheses), and then (2) says that we can get any free filter of $\ssy(\MM)$. We will be interested in getting free ultrafilters. One of the problems that occurs is how to distinguish the V-point $a$ that we designate as the generator of $\NN$ over $\MM$ from other V-points that generate $\NN$ over $\MM$. Or, more precisely, how do we distinguish the (ultra)filter ${\mathcal F}_a$ from other possible filters.
  The next lemma will be a help in doing this. 
  
  \bigskip

{\sc Lemma 4.5}: {\em Suppose that $Y$ is $I$-viable,  $d$ is nonstandard, and $\langle f_i : i < d \rangle$ is an $\MM$-definable $d$-tuple of functions on $M$.  Then there are nonstandard $j \leq d$ and an $(I \cap [0,j))$-viable  $Z \subseteq Y$  such that if $i \in [0,j)$, then $f_i$  is canonical on~$Z$.}

   \bigskip
 
 {\it Proof}. Without loss of generality, assume that $I = [0,k)$ and $k = d$. Define $C^j_k(v)$ for standard $j$ so that $C^0_k(v) = v$ and $C^{j+1}_k(v) = C(C^j_k(v))$. Since $J$ is exponentially closed, we can also assume that $k$ is small enough so that there is $v \in M$ such that
 $$
 J < v < C^k_k(v) < v^v < w.
 $$
 For each $j$, where $\min(I) < j < \omega$, there is, by repeated applications of Corollary~3.7(\pa), an $( [0,j),v)$-special $Z_j \subseteq Y$ such that for each $i \leq j$, $f_i$ is  canonical on $Z_j$. By overspill, there is a nonstandard $j \leq k$ and an $([0,j),v)$-viable $Z \subseteq Y$ such that for each $i \in  [0,j)$, $f_i$ is canonical on $Z$. \qed

    \bigskip
 
 Having the previous corollary, we  now add a fifth condition to  (S1) -- (S4). Keep in mind that these conditions, as stated in the proofs of Theorems~1.1 and~1.2, should be applied to a sequence of viable (and not just large) sets. 
 
 \smallskip
 
 \begin{itemize}
 
 \item[(S5)] If   $d$ is nonstandard and $\langle f_j  : j < d \rangle$ is an $\MM$-definable $d$-tuple of functions  on $M$, then there are $i < \omega$ and nonstandard $k \leq d$ such that for every $j < k$, $f_j$ is canonical on $X_i$. 
 
 \end{itemize}
 
 \smallskip
 \noindent We will say that an extension $\NN$ of $\MM$ is a {\bf V$^*$-{extension}} if, in addition to being a V-extension, it also incorporates (S5). We also define {\bf V$^*$-point} in the obvious way.

    \bigskip
 
 {\sc Lemma 4.6}: {\em Suppose that $\NN$ is a V$^*$-extension of $\MM$ and that $a,b \in N$ are V$^*$-points.   Then ${\mathcal F}_b \cong_{\sf RK}^\MM {\mathcal F}_a$.}

 \bigskip
 
  {\it Proof}. By symmetry, it suffices to show that ${\mathcal F}_b \leq_{\sf RK}^\MM {\mathcal F}_a$. 
  
  Let $a \in P_h^\NN$ and $b \in P^\NN_{h'}$, where $(P_h; \overline E)$ and $(P_{h'};\overline {E'})$ are basic $h$- and $h'$-structures, respectively.  Let $X_0 \supseteq X_1 \supseteq X_2 \supseteq \cdots$ be a sequence demonstrating that $a$ is V$^*$-point, and $Y_0 \supseteq Y_1 \supseteq Y_2 \supseteq \cdots$ be a sequence demonstrating that $b$ is V$^*$-point. Since both $a$ and $b$ generate $\NN$ over $\MM$, we can let $f : M \into M$ be an $\MM$-definable, one-to-one function such that $f^\NN(a) = b$.  It's safe to assume that $X_0 \subseteq f^{-1}(Y_0)$. 
  
Suppose that $i_0,i_1 < \omega$ are such that $f : X_{i_0} \into Y_{i_1}$. For each $j \leq h'$, let 
 $h_{i_0,i_1,j} : X_{i_0} \into M$ be such that if $x \in X_{i_0}$, then $h_{i_0,i_1,j}(x)$ is the equivalence class of $E_j \cap Y^2_{i_1}$ to which $f(x)$ belongs.   By (S5), we get $i < \omega$ and    nonstandard $k \leq h'$    such that $i \geq i_0$ and for each $j < k$, $h_{i_0,i_1,j}$ is canonical on $X_i$. Thus, there are $d_0 \leq d_1 \leq d_2 \leq \cdots \leq d_{k-1}$, where  each $d_j$ is maximal such that $h_{i_0,i_1j}$ is $d_j$-canonical on $X_i$. Hence, for all $x,y \in X_i$ and $j < k$,
  $$
  \langle x,y \rangle \in E_{d_j} \Longleftrightarrow \langle f(x),f(y) \rangle \in E'_j.
  $$
   We claim: $\{d_j : j < \omega\}$ is an infinite subset of $\omega$.
   
   For notational simplicity, we prove the claim only in the case that $i_0 = i_1 = 0$, and then we let $h_j = h_{0,0,j}$.

   First, we show that $d_j < \omega$ for each $j < \omega$.  For a contradiction, suppose that  $j < \omega$ and $d_j$ is nonstandard.  Let $i_1 < \omega$ be such that $Y_{i_1}$ is $I_1$-viable and $\min(I_1) > j$. Let $i_0 < \omega$   be such that $i_0 \geq i$, $X_{i_0} \subseteq f^{-1}(Y_{i_1})$ and $X_{i_0}$ is $I_0$-viable, where $\max(I_0) < d_j$. Then, $E_{d_j} \cap X_{i_0}^2$ is discrete and $E'_j \cap Y_{i_1}^2$ is trivial.  Thus, 
  \begin{align*}
  x,y \in X_{i_0} & \Longrightarrow  f(x),f(y) \in Y_{i_1} \\
                         & \Longrightarrow  \langle f(x),f(y) \rangle \in E'_j \\
                         & \Longrightarrow \langle x,y \rangle \in E_{d_j} \\
                         & \Longrightarrow x = y,
  \end{align*}
  contradicting that $X_{i_0}$ is infinite. 
  
   Next, we show that $\{d_j : j < \omega\}$ is infinite. For a contradiction, suppose that $d < \omega$ is such that $d_j < d$ for all $j < \omega$. By overspill, we can assume that $d_j < d$ for all $j < k$. Let $i_0 < \omega$ be such that $i_0 \geq i$, $X_{i_0}$ is $I_0$-viable, $\min(I_0) > d$ and $X_{i_0} \subseteq f^{-1}(P_{h'})$. Let $i_1$ be such that $Y_{i_1} \subseteq f[X_{i_0}]$, and let $Y_{i_1}$ be $I_1$-viable. Then, 
    \begin{align*}
  x,y \in X_{i_0} & \Longrightarrow  \langle x,y \rangle \in E_d \\
 & \Longrightarrow \langle x,y \rangle \in E_{d_j} {\mbox{ for all }} j < \omega \\ 
                         & \Longrightarrow  \langle f(x),f(y) \rangle \in E'_j {\mbox{ for all }} j < \omega\\
                         & \Longrightarrow \min({I_1}) {\mbox{ is nonstandard,}}
  \end{align*}
  contradicting that $I_1 \cap \omega$ is infinite. The claim is proved.
 
 Let $i < \omega$ and $k \leq h'$ be as in the claim. Let $D = \{d_j : j < \omega\}$. We have just proved that $D$ is an infinite subset of $\omega$. If  $X_i$ is $I$-viable, then $I \cap \omega \subseteq D$, so that $D \in {\mathcal F}_a$.
  
  Let $g : \omega \into \omega$ be such that $g(n) = j$ iff $j$  is the least such that if $ n \leq d_j$. We have just seen that $g$ is well defined. Obviously, $g \in \ssy(\MM)$. 
 We next claim that $g$ demonstrates that ${\mathcal F}_b \leq_{\sf RK}^\MM {\mathcal F}_a$.
 
 Suppose not. Let $B \in {\mathcal F}_b$ be such that $g^{-1}(B) \not\in {\mathcal F}_a$. Let $i_1 < \omega$ be such that $Y_{i_1}$ is $I_1$-viable and $B \supseteq I_1 \cap \omega$. 
 Since $g^{-1}(B) \in \ssy(\MM)$ and ${\mathcal F}_a$ is an ultrafilter of $\ssy(\MM)$, it must that 
 $\omega \backslash g^{-1}(B   ) \in {\mathcal F}_a$. Let $i_0 < \omega$ be  such that  $i_0 \geq i$, 
 $X_{i_0} \subseteq f^{-1}(Y_{i_1}$, $X_{i_0}$ is $I_0$-viable and  $I_0 \cap g^{-1}(B) = \varnothing$. Then $h_{i_0,i_1,0}(0)$ is nonstandard, which contradicts the   claim.  \qed
 
 \bigskip

We complete the proof of Theorem~B. Let $I \subseteq \omega$. Let $\lambda_I$ be as in the proof of Corollary~2.12, and then, by Theorem~2.11,  let $\NN_0$ be a humble extension of $\MM$, generated by $c$ over $\MM$,  such that  
$\Lt_0(\NN_0 / \MM) \cong 2 + \lambda_I + 1$. Let ${\mathcal U}$ be an ultrafilter of $\ssy(\MM)$ and, by  Proposition~2.3(2), $\MM^\omega/{\mathcal U} \cong \NN_0$ with an isomorphism that fixes $M$.
Let $\NN$ be a V$^*$-extension of $\MM$ generated over $\MM$ by the V$^*$-point $a$ such that 
${\mathcal F}_a = {\mathcal U}$. If $b \in N$ is any V$^*$-point that generates $\NN$ over $\MM$, then 
${\mathcal F}_a \cong_{\sf RK}^\MM {\mathcal F}_b$ by Lemma~4.6. Thus, $\NN_0 \cong \MM^\omega/{\mathcal F}_a \cong \MM^\omega /{\mathcal F}_b$, and the isomorphisms fix $M$. Since $\MM^\omega /{\mathcal F}_b$ is definable in $(\NN,M,b)$, we get, from Corollary~2.12, that $I$ is definable in $(\NN,M)$.
Clearly, there is a formula $\varphi(x)$ that works uniformly for all $\MM$, $J$ and $I$, completing 
the proof of Theorem~B. 

%
%

\bibliographystyle{plain}

\end{document}